\documentclass[11 pt]{article}

\usepackage[latin1]{inputenc}
\usepackage{amsmath,indentfirst,qsymbols,graphicx,psfrag,amsfonts,stmaryrd,hyperref}

\DeclareMathOperator{\Id}{Id}

\DeclareMathOperator{\Trace}{Trace}

\DeclareMathOperator{\Div}{div}
\DeclareMathOperator{\Mod}{mod}

\def \1{\textbf{1}}

\def \Ch{{\sf Ch}}

\def \bB{{\bf B}}

\def \bW{{\bf W}}
\def \bw{{\bf w}}

\def \bQ{{\bf Q}}
\def \bq{{\bf q}}
\def \bP{{\bf P}}
\def \bX{{\bf X}}
\def \bY{{\bf Y}}
\def \bZ{{\bf Z}}

\def \bT{{\bf T}}
\def \bM{{\bf M}}

\def \bW{{\bf W}}
\def \bar{\overline}
\def \ben{\begin{eqnarray}}
\def \een{\end{eqnarray}}
\def \ba{\begin{align}}
\def \ea{\end{align}}
\def \matz{\begin{bmatrix}0\end{bmatrix}}

\def \be{\begin{eqnarray*}}
\def \ee{\end{eqnarray*}}
\def \bh{{\bf h}}

\def \beq{\begin{equation}}
\def \eq{\end{equation}}

\def \build#1#2#3{\mathrel{\mathop{\kern 0pt#1}\limits_{#2}^{#3}}}

\def \ba{{\bf a}}
\def \bb{{\bf b}}
\def \bv{{\bf v}}
\def \bx{{\bf x}}
\def \by{{\bf y}}
\def \bz{{\bf z}}

\def \captionn#1{\begin{center}\begin{minipage}{15cm}\sf\caption{\small #1}\end{minipage}\end{center}}
\def \Sq{{\sf Sq}}

\def \Sys{{\sf Sys }}

\def \Proof{{\bf Proof. \rm}}
\def \Proofof#1{{\bf Proof of #1. \rm}}

\def \eref#1{(\ref{#1})}

\def \G{{\bf GF}}
\def \A{{\bf A}}
\def \imp{\Rightarrow}

\def \l{\left}

\def \r{\right}
\def \sous#1#2{\mathrel{\mathop{\kern 0pt#1}\limits_{#2}}}
\def \sur#1#2{\mathrel{\mathop{\kern 0pt#1}\limits^{#2}}}

\def \row {{\sf Row }}
\def \nn{_{n}}
\def \tb{{\tiny $\blacksquare$~ }}

\newqsymbol{`B}{\mathcal{B}}
\newqsymbol{`E}{\mathbb{E}}
\newqsymbol{`I}{\mathbb{I}}
\newqsymbol{`N}{\mathbb{N}}
\newqsymbol{`O}{\Omega}
\newqsymbol{`P}{\mathbb{P}}
\newqsymbol{`Q}{\mathbb{Q}}
\newqsymbol{`R}{\mathbb{R}}
\newqsymbol{`W}{\mathbb{W}}
\newqsymbol{`Z}{\mathbb{Z}}
\newqsymbol{`a}{\alpha}
\newqsymbol{`e}{\varepsilon}
\newqsymbol{`o}{\omega}
\newqsymbol{`t}{\tau}
\newqsymbol{`w}{{\cal W}}
\def\cro#1{\llbracket#1\rrbracket}

\begin{document}

\newtheorem{fig}{\hspace{2cm} Figure}
\newtheorem{lem}{Lemma}
\newtheorem{defi}[lem]{Definition}
\newtheorem{pro}[lem]{Proposition}
\newtheorem{theo}[lem]{Theorem}
\newtheorem{cor}[lem]{Corollary}
\newtheorem{note}[lem]{Note}
\newtheorem{conj}{Conjecture}
\newtheorem{Ques}{Question}

\renewcommand{\baselinestretch}{1.05}

\begin{center}
\huge\bf
Directed animals, quadratic and rewriting systems\\
{\large \bf Jean-Fran\c{c}ois Marckert}
\rm \\
\large{CNRS, LaBRI, Universit\'e Bordeaux \\
 351 cours de la Libération\\
33405 Talence cedex, France}
\normalsize
\end{center}

\begin{abstract} A directed animal is a percolation cluster in the directed site percolation model. The aim of this paper is to exhibit a strong relation between the problem of computing the generating function $\G$ of directed animals on the square lattice, counted according to the area and the perimeter, and the problem of solving a system of quadratic equations involving unknown matrices. We present some solid evidence that some infinite explicit matrices, the fixed points of a rewriting like system are the natural solutions to this system of equations: some strong evidence is given that the problem of finding $\G$ reduces to the problem of finding an eigenvector to an explicit infinite matrix. 
Similar properties are shown for other combinatorial questions concerning directed animals, and for different lattices. 
\end{abstract}

\begin{quotation}{\footnotesize 
The author is partially supported by the ANR-08-BLAN-0190-04 A3. }
\end{quotation}
%\tableofcontents

\section{Introduction}
\normalsize

We are mainly interested in the study of the area and perimeter generating function $\G^\Sq$ of directed animals on the square lattice $\Sq$, but other lattices and questions will also be addressed. The computation of  $\G^\Sq$ is a central question in enumeration problems for directed animals on two dimensional lattices, since it is deeply related to the study of directed percolation on the square lattice. In this paper, even if we do not find an explicit formula for $\G^\Sq$, we show that to compute $\G^\Sq$ it suffices to solve a quadratic system of equations involving 4 unknown finite matrices. We are unable to find a solution, but we provide some infinite size matrices which appear as the natural solution to this system of equations. They appear to be a fixed point of a rewriting system, the rewriting rules involving the tensorial product of matrices. We give strong evidence that finding a right and a left eigenvector to these matrices should lead to $\G^\Sq$.  We hope that this gives some insight on the algebraic structure of this problem, and that this will allow some readers to compute $\G^\Sq$. \par
In Section \ref{sec:extensions}, we show that numerous similar problems can be treated similarly. \medskip
\begin{figure}[ht]
\psfrag{l0}{$l_0$}
\psfrag{l1}{$l_1$}
\psfrag{l2}{$l_2$}
\centerline{\includegraphics[height=3cm]{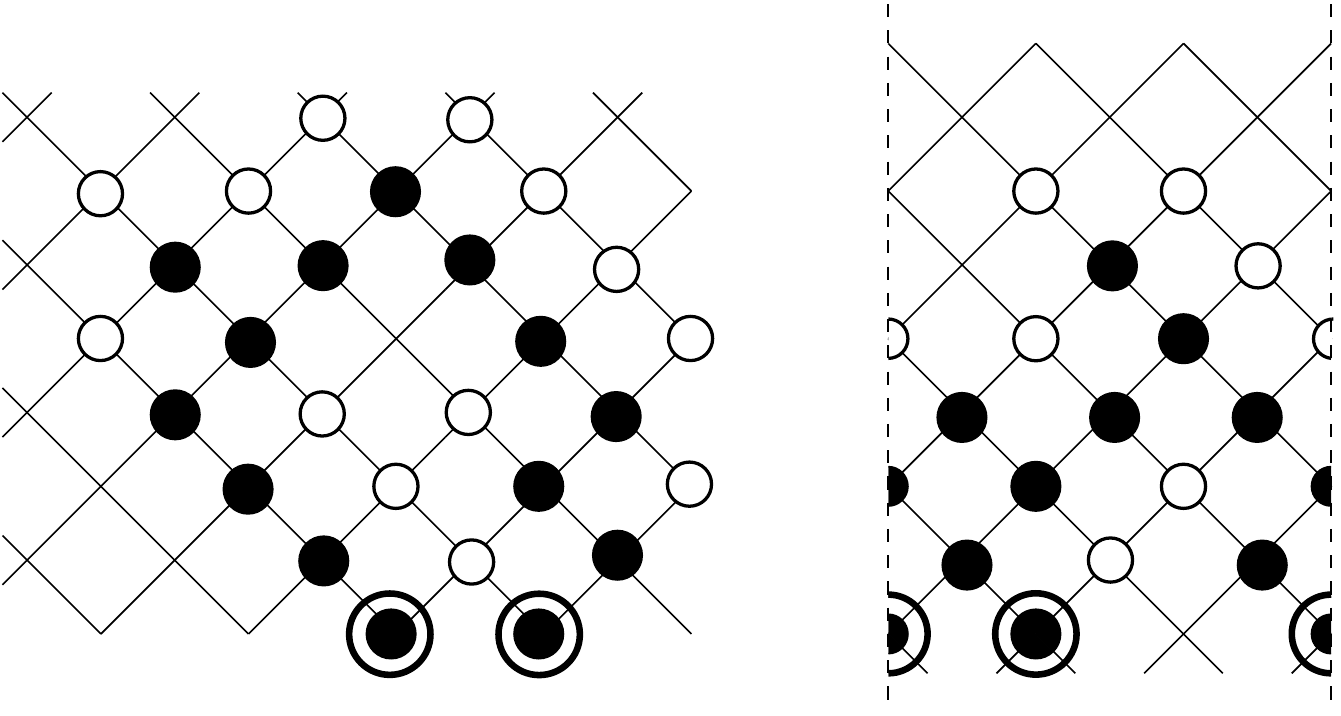}}
\captionn{\label{fig:STH}DA on the square lattice and on the cylinder $\Sq(3)$. The black dots are the cells and  the white dots, the perimeter sites. The cylinder is defined as the square lattice with periodic conditions, meaning that the extreme vertical lines on the second picture are identified. The outlined cells constitute the source.}
\end{figure}

The set of oriented graphs with no cycle and no multiple edges which have a finite or countable number of vertices and bounded degree is denoted ${\cal G}$. For any graph $G=(V,E)$ in ${\cal G}$, $V$ is the set of vertices and $E\subset V^2$ the set of oriented edges. 
The orientation of the edges leads to the notion of a descendant: for $(x_1,x_2) \in E$, $x_2$ is said to be a child of $x_1$ and the set of children of $x_1$ is denoted $\Ch(x_1)$. A directed path $d$ in $G$ is a sequence of vertices $(x_1,\dots,x_k)$ such that for any $l\geq 2$, $x_l\in \Ch(x_{l-1})$. The vertex $x_1$ (resp. $x_k$) is called the origin (resp. the target) of $d$. 
\begin{defi}\label{diran} Let $G=(V,E)$ be in ${\cal G}$, and $S$ be a subset of $V$. \\
$\bullet$ A directed animal (DA) $A$ with source $S$ is a subset of $V$ containing $S$, such that for every $a\in A$ there exists a directed path having target $a$ and its origin in $S$ entirely contained in $A$. The cardinality $\#A$ of $A$ is called the area of $A$.\\
$\bullet$ A perimeter site $c$ of a DA $A$ with source $S$ is an element of $V\setminus A$ such that $\{c\} \cup A$ is still a DA with source $S$. The set of perimeter sites of $A$ is denoted $P(A)$.\\
\end{defi}
We denote by ${\cal A}^G_S$ the set of finite DA on $G$ with source $S$. The generating function (GF) $\G_S^G$ counts the DA with source $S$ according to the area and perimeter:
\[\G_S^G(x,y):=\sum_{A\in {\cal A}^G_S}x^{\# A}y^{\#P(A)}.\]
Hence, the area generating function is $\G_S^G(x,1)$. \par
The search for a formula for $\G_S^G(x,y)$ may be seen as the combinatorial contribution to the study of directed percolation per site models. Indeed, on a probability space $(\Omega,{\cal A},`P)$ consider a random colouring of the vertices of $V$ by the colours 0 and 1. Formally, this is given by a family of i.i.d. Bernoulli random variables $(B^v(p),v\in V)$ indexed by the vertex set (we then have $`P(B^v(p)=1)=1-`P(B^v(p)=0)=p$).  The directed percolation cluster with source $v \in V$ is the maximum DA with source $S=\{v\}$ included in the set of 1-coloured vertices, that is $\{u \in V~:B^u(p)=1\}$ (the empty case, possible here, arises with probability $1-p$): denote it ${\bf A}^v(p)$. Since for any DA A with source $v$,
$`P({\bf A}^v(p)=A)=p^{\#A}(1-p)^{\#P(A)}$, the percolation cluster is finite with probability 1 if  $(1-p)+\sum_{A\in{\cal A}^G_{\{v\}}} p^{\#A}(1-p)^{\#P(A)}=1$, which is equivalent to $\G_{\{v\}}^G(p,1-p)=p$. Hence a computation of $\G_{\{v\}}^G(x,y)$ would probably allows one to compute the directed percolation threshold, and/or the associated critical exponent. \medskip

Denote by $\Sq=(V_\Sq,E_\Sq)$ the directed square lattice where $V_\Sq=\mathbb{Z}^2$ and \[E_\Sq=\l\{((x_1,y_1),(x_2,y_2))\in V_\Sq^2, \textrm{ such that } (x_2,y_2)-(x_1,y_1)\in\{(0,1),(1,0)\}\r\}.\]
Surveys on the study of DA  on two dimensional lattices exist:  Bousquet-Mélou \cite{BM1} and Le Borgne \& Marckert \cite{LB-M}.  When $S$ is reduced to a singleton, the area GF $\G_S^{\Sq}(x,1)$ is well known, and numerous different approaches are possible to compute it: the gas approach (Dhar \cite{DH1,DH5}, but also \cite{BM1}, \cite{LB-M}, Albenque \cite{Al}), heap of pieces approach (Viennot \cite{VI1}), combinatorial decomposition (Corteel \& al \cite{CDG}, Bétréma \& Penaud \cite{BEPE}). On the other side, almost nothing is known about $\G_S^{\Sq}(x,y)$ (except for Bacher  \cite{BACH1} who computed $\sum p(A)x^{|A|}=\frac{\partial \G_S^{\Sq}(x,u)}{\partial u}(x,1)$ on the square lattice with or without periodic conditions, proving conjectures by Conway \cite{CC} and Le Borgne \cite{LB2}). $\G_S^{\Sq}(x,y)$ is not believed to be $D$-finite.\medskip

The aim of this paper is to use some algebra to search for a formula for $\G_S^{\Sq}(x,y)$. 
We will use the idea of  Nadal \& al. \cite{NDV} and Hakim and Nadal \cite{HV}, also  used extensively  in \cite{BM1}. First, the work is done on a so-called cylinder $\Sq(n)$, a vertical strip of $\Sq$ with periodic conditions (see Figure \ref{fig:STH}). Second, the corresponding for $\Sq$ is obtained by taking a formal limit since small DA on $\Sq(n)$ and $\Sq$ are the same. We will proceed similarly here, starting with $\G^{\Sq(n)}(x,y)$.\par
In order to highlight the different considerations leading to the introduction of infinite matrices, we have decided to simultaneously treat the study of $\G^{\Sq(n)}(x,y)$ and a case where infinite matrices can be avoided: the computation of $\G^{\Sq(n)}(x,1)$. This leads to a new derivation of $\G^{\Sq}(x,1)$.  

\subsection{Two gases}
Let  $G=(V,E)\in{\cal G}$ be an oriented graph. Following the ideas developed in \cite{LB-M}, we define two processes $\bX_G=(X_v,v \in V)$ and $\bY_G=(Y_v, v\in V)$ indexed by the vertex set, and taking their values in $\{0,1\}$. For this latter reason, the processes are called ``gases'', the value 1 (resp. 0) representing the presence (resp. absence) of a particle.\par
Both processes $\bX_G$ and $\bY_G$ are defined on a probability space $(\Omega,{\cal A},`P)$, on which are defined some families of i.i.d. random variables $(\Xi^v,v\in V)$ indexed by the vertex set, where $\Xi^v=(B_1^v(p),B_2^v(q))$ is a pair of independent Bernoulli random variables whose parameters are $p$ and $q$.
 
\paragraph{Gas of type 1 :} For any $v\in V$, set
\begin{equation}\label{eq:gas1}
X_v=B^v_1(p)\prod_{c \in \Ch(v)} (1-X_c). 
\end{equation}
That is, if $X_c=0$ for all $c\in \Ch(v)$, then $X_v=1$ with probability $p$;  otherwise $X_v=0$. Since two neighbouring sites can not be simultaneously occupied, this model is called a hard particle model in the physics literature. 

\paragraph{Gas of type 2 :} For any $v\in V$, set
\begin{equation}\label{eq:gas2}
Y_v=B^v_1(p)\min\{Y_c~: c \in \Ch(v)\} +(1-B^v_1(p))B_2^v(q). 
\end{equation}
Here $Y_v$ is equal to $\min\{Y_c~: c\in \Ch(v)\}$ with probability $p$, and to  $B_2^v(q)$. with probability $1-p$.

\begin{lem}\label{lem:welldefi}
Let $G\in {\cal G}$. If $p\in(0,1)$ is small enough, both processes $\bX_G$ and $\bY_G$ are almost surely well defined.
\end{lem}
\Proof We use the argument in \cite{LB-M} (the argument being already present in the PhD thesis of Le Borgne \cite{LBPhD}). For both gases, when $B^v_1(p)=1$, the set of values $\{X_c,c\in \Ch(v)\}$ (resp. $\{Y_c,c\in \Ch(v)\}$) is needed to compute $X_v$ (resp. $Y_v$), but they are not needed when $B^v_1(p)=0$, in which case $X_v=0$ and $Y_v=B_2^v(q)$. The fact that these ``recursive definitions'' \eref{eq:gas1} and \eref{eq:gas2} indeed define some objects is not clear, but the values of both $X_v$ and $Y_v$ are certainly well defined if ${\bf A}^v(p)$ is finite, since in this case, the recursive computation of $X_v$ (and $Y_v$) using the values of the children ends since the value of $X$ and $Y$ on perimeter sites of  ${\bf A}^v(p)$  -- sites where $B_1(p)$ is zero -- is well defined. Hence, if the family of DA $({\bf A}^v(p),v\in V)$ is a family of finite DA, both processes are defined. Now consider the standard problem of the directed percolation threshold. Let
\[p_{crit}=\sup\{p~: `P(\forall v, |{\bf A}^v(p)|<+\infty)=1\}.\]
Since we assume that the maximum degree if the graph is finite, $p_{crit}\in(0,1]$. For all $0\leq p<p_{crit}$, $\bX_G$ and $\bY_G$ are then a.s. defined. ~$\Box$ \medskip 

A  subset $S$ of $V$ is said to be free if for any $s_1,s_2 \in S$ with $s_1\neq s_2$, there does not exist any directed path with origin $s_1$ and target $s_2$ in $G$. The following Proposition says that the computation of the finite dimensional distribution of the gas of type 1 (resp. type 2) is equivalent to the computation of the DA GF according to the area (resp. area and perimeter), for general source. A DA is said to have over-source $S$ if is a DA with source $S'$, with $S'\subset S$. If $S'\neq S$, the set $S\setminus S'$ is taken to be a subset of the perimeter. We denote by $\bar{\G_S^G}$ the generating function of DA with over-source $S$.
\begin{pro}\label{pro:fo} 
For any directed graph $G=(V,E)$ in ${\cal G}$ and any free subset $S$ of $V$,  
\be
`P(X_v=1,v\in S)&=&(-1)^{\#S}\G_S^G(-p,1),\\
`P(Y_v=1,v\in S)&=&\bar{\G_S^G}(p,(1-p)q),
\ee
where the first equality holds for $|p|$ smaller than the radius of convergence of $\G_S^G(-p,1)$ and the second one holds if $0 < p< p_{crit}$.  
\end{pro}
\Proof The first assertion is proved in \cite{LB-M} (Theorem 2.7 for a single source, Proposition 2.16 for any source) on a general graph, and was already used in  \cite{BM1} on lattices for a single source. Dhar \cite{DH1}, who made the connection on lattices between GF of DA and the problem of finding the density of a hard particle system on an associated graph, did not use the construction of the process $X$, but different considerations of the same process. The second assertion is well-known, and it is also proved in \cite{LB-M} (Theorem 4.3) and valid on any graph of ${\cal G}$. The reason is simple: $Y_v=1$ if and only if $B^u_2(q)=1$ for all perimeter sites $u$ of $\A^v(p)$ (where $v$ is considered as a perimeter site of $\A^v(p)$ in the case where $\A^v(p)$ is empty). ~$\Box$ 

\begin{note} Although $\G(x,1)$ is a projection of $\G(x,y)$, and $\bar{\G_S^G}(x,y)$ can be computed easily thanks to 
$(\G_{S'}^G(x,y),S'\subset S)$,  the gas of type 2 is not an extension of the gas of type 1. To compute $\G^G(p,1)$ using the gas of type 2, $q$ needs to be $1/(1-p)$, which is larger than 1; this is not possible for probabilistic considerations. Nevertheless, given the polynomial form of $\bT^{\bY}$ (see Formulas \eref{main:equa1} and \eref{eq:trans}), $\bT^\bY$ still has a meaning when $q=1/(1-p)$; for this value, it is no longer a positive kernel, but $\sum_{\bb \in E_n} \bT^{\bY}_{\ba,\bb}=1$ for any $\ba\in E_n$. A non-negative solution $\mu$ to $\mu=\mu\bT^{\bY}$ still exists, as can be checked by rewriting the following system of equations: let, for any $C \subset\{0,\dots,n-1\}$, 
\[W_C=\mu(\{x \in E_n ~:   i\in C \imp x_i= 1\}).\]
$\mu$ solves the system $\mu=\mu\bT^{\bY}$ if and only if $W:=(W_C,C \subset\{0,\dots,n-1\})$ is solution to
$W_C =\sum_{D\subset C} ((1-p)q)^{C\setminus D} p^D W_{\Ch(D)}$,
which when $(1-p)q=1$ is a rewriting of $W_C =\sum_{D\subset C}  p^D W_{\Ch(D)}$. Clearly, $W$ satisfies the same system of 
equations as $\bar{\G_C^G}:=(\bar{\G_C^G}(p,1),C \subset\{0,\dots,n-1\})$. Since $\bar{\G_C^G}$ exists and is non negative, $\mu=\mu\bT^{\bY}$ admits some solutions. The initial conditions $W_\varnothing=\bar{\G_\varnothing^G}(p,1)=1$ allows one to identify the two set of series $W$ and $\bar{\G_C^G}$.
As in Dhar \cite{DH1} or Bousquet-Mélou \cite{BM1}, $W_C$ has a product form on $E_n$, meaning that
 ($W_\ba=\alpha_n \prod_{i=1}^n Q_{a_i,a_{i+1\mod n}}$ for some numbers $Q_{0,0},Q_{0,1},Q_{1,0},Q_{1,1}$). The behaviour of $`P(Y_v=1,v\in S)$ becomes singular at $q=1/(1-p)$:  some drastic simplifications of the involved algebra appear, but only for that value of $q$. This leads  to a Markovian type structure of $\mu$. \par
Some extensions of the gas of type 1 can also be introduced, for example 
$X_x=B_p\prod_{c\in\Ch(x)}(1-X_c)+(1-B_p)B_q$. These extensions are either of the same type as this gas, meaning that an easy-to-prove Markovian behaviour occurs, or will present  the same kind of difficulties as for the gas of type 2 (as illustrated by the different cases discussed in Section \ref{sec:extensions}). We think that any gas allowing us to compute $\G^{\Sq(n)}(x,y)$ will be (at best) as difficult to describe as the gas of type 2.
\end{note}

\subsection{The square lattice : the cylinder approach}

We study some properties of the processes of type 1 and 2 defined on $\Sq(n)$, the square lattice with periodic conditions shown in Figure \ref{fig:STH}. We first need to label the sites of this lattice and its rows, in order to describe its Markov chain structure, row by row: the  $l$th row 
is $\row^{(l)}=\l\{(x \mod n,l-x \mod n),x \in \mathbb{Z}\r\}$ and let $x^{(l)}(i)=(i \mod n,l-i\mod n)$ be the $i$th element of this row. The two children of $x^{(l)}(i)$ are  $x^{(l+1)}(i+1)$ and $x^{(l+1)}(i)$. Note that $\row^{(l+1)}$ is above $\row^{(l)}$ in Figure \ref{fig:STH}.\par
It is useful to provide a row decomposition of the processes $\bX_{\Sq(n)}$ and $\bY_{\Sq(n)}$. For brievity, let $\bZ\in\{\bX_{\Sq(n)},\bY_{\Sq(n)}\}$ be one of these processes, and let $Z^{(l)}(i)$ be the value of $\bZ$ at $x^{(l)}(i)$, and
\[\bZ^{(l)}=\l(Z^{(l)}(i),i=0,\dots,n-1\r)\] 
the value of $\bZ$ on $\row^{(l)}$. Clearly $\bZ^{(l)}$ depends on  $\bZ^{(l+1)}$ and the values of the Bernoulli random variables $(\Xi_v,v\in \row^{(l)})$ only.  Then, the sequence $\bZ^{(l)}$, when $l$ goes from $+\infty$ to $-\infty$, is a Markov chain. The transition kernel $\bT_n^\bZ$ can be expressed using  \eref{eq:gas1} and \eref{eq:gas2}~: denote by $E_n$ the state space  $\{0,1\}^n$, by $\ba\nn:=(a_0,\dots,a_{n-1})$, and $\bb\nn=(b_0,\dots,b_{n-1})$ some generic elements of $E_n$, and by $\oplus$ the addition in $\mathbb{Z}/n\mathbb{Z}$. We have
\begin{equation}\label{main:equa1}\bT^\bZ(\bb,\ba):=`P(\bZ^{(l)}=\ba ~|~\bZ^{(l+1)}=\bb)=\prod_{i=0}^{n-1} T^{\bZ}_{b_i,b_{i\oplus1},a_{i}},
\end{equation}
and where $T^{\bX}_{x,y,z}$ and $T^{\bY}_{x,y,z}$ are explicit:
\ben\label{eq:trans}
T^{\bX}_{x,y,z}&=& `P(B_1(p)(1-x)(1-y)=z)\\
               &=&1_{(x,y)\neq (0,0)}1_{z=0}+1_{(x,y)=0}\l( p 1_{z=1}+(1-p)1_{z=0}\r),\\ 
T^{\bY}_{x,y,z}&=& `P(B_1(p)xy+(1-B_1(p))B_2(q)=z)\\
\nonumber&=&1_{\min(x,y)=0}\l( (p+(1-p)(1-q))1_{z=0}+(1-p)q1_{z=1}\r)\\
&&+1_{\min(x,y)=1}\l( (p+(1-p)q)1_{z=1}+(1-p)(1-q)1_{z=0}\r). 
\een
If, for a given $l\in\mathbb{Z}$, $\bZ^{(l+1)}$ has distribution $\nu$, then the distribution $\mu$ of $\bZ^{(l)}$ satisfies 
\begin{equation}\label{main:equa0}
\mu(\ba)=\sum_{\bb\in E_n} \nu(\bb)\,  \bT^{\bZ}_n(\bb,\ba).
\end{equation}
For clarity, we use an operator-type notation and write $\mu=\nu\,\bT^\bZ_n$; accordingly,  $\mu=\nu\,(\bT^\bZ_n)^\kappa$ designates the law of $\bZ^{(-\kappa)}$ if $\bZ^{(0)}$ has distribution $\nu$. \par
Using the row decomposition of $\Sq(n)$, we have 
\begin{lem}\label{lem:MC} 1) Let $n\geq 1$ be fixed and $\bZ\in\{\bX_{\Sq(n)},\bY_{\Sq(n)}\}$.  For $p,q\in (0,1)$,  the process $(\bZ^ {(-j)},j\in \mathbb{Z})$ (indexed by decreasing $j$'s) is a Markov chain with finite state spaces $E_n$ under its stationary distribution. Since this Markov chain is irreducible and aperiodic, the distribution of $\bZ^{(j)}$ for any $j$ is the only probability measure solution $\mu^\bZ$ of
\begin{equation}\label{main:equa}
\mu^\bZ=\mu^\bZ\bT^\bZ_n.
\end{equation}
2) The distribution $\mu^\bZ$ is the limit distribution of any Markov chain with transition $\bT^\bZ_n$ on $\Sq(n)$, i.e. for any distribution $\nu$ on $E_n$, we have  $\nu (\bT_n^\bZ)^\kappa\sous{\longrightarrow}{\kappa\to +\infty}\mu^\bZ$ .
\end{lem}
The matrix $\bT^Z_n$, indexed by the elements of $E_n$ is called a (row) transfer matrix in statistical physics literature.  By the Perron-Frobenius theorem, the equation \eref{main:equa} has a unique solution up to a multiplicative constant, which can be turned into a probability distribution. Finding a $\mu^\bZ$ solution of \eref{main:equa} is a problem of linear algebra, and a solution can be found by using a computer for small values of $n$. The fact that $\bT_n^\bZ$ has the product representation given in \eref{main:equa1} plays a secondary role in that respect. In the following, it will play a role of primary importance. \medskip

\noindent\Proof The only unclear assertion in (1) is that the row process $\bZ^{(-l)}$ follows under its stationary distribution. This comes from the infinite construction that we have, under which $\bZ\nn^{(l)}$ and $\bZ\nn^{(l+1)}$ have the same distribution. Assertion (2) is a well known property for aperiodic irreducible Markov chains on a finite state space.~$\Box$ 

\begin{note}
In the case of the square (or triangular) lattice the gas of type 1 is a Markov chain  on the horizontal lines of the lattice \cite{LB-M}. In \cite{BM1}, it is observed that it is a Markovian field on a zigzag formed with two consecutive lines of $\Sq(n)$; this Markovian field  converges when the row size goes to $+\infty$ to a Markov chain on the line according to \cite{Al}. 
It turns out that the gas of type 2 is not Markovian on a line of the cylinder, on a line of the lattice, or on the zig-zag \cite{BM1}. In other words, there does not exist any $2\times 2$ matrix $A$ such that
\beq\label{eq:prodd}
\mu^{\bY}_{n}(x_1,\dots,x_n)=  c_n \prod_{i=1}^n A_{x_i,x_{i\otimes +1}}.
\eq
To show this, one can solve \eref{main:equa} with a computer for $n=3,4,\dots$, to check that no factorisation of the invariant measure compatible with \eref{eq:prodd} is possible (a similar work can be done directly on the rows of the entire lattice). It can also be checked that no solution corresponding to a Markov chain with longer memory exists, at least for small memories. When a computer is used to compute the invariant distribution on a cylinder of small size, no ``regularity'' of the measure has been observed: it turns out that on the cylinder of size 12, if $\mu_{12}^{\bY}(x_1,\dots,x_n)=\mu_{12}^{\bY}(y_1,\dots,y_n)$, then $\bx$ and $\by$ are equal up to a rotation and/or symmetry. This implies that what happens is drastically more complex that what appears for the gas of type 1 where $\mu^\bX_n$ has the form \eref{eq:prodd}.
\end{note}
We now look at new considerations.

\section{A new paradigm}

Proposition \ref{pro:fo} says that  $\mu_n^\bY(1,\star,\dots,\star):=\sum_{x_i\in\{0,1\},i=2,\dots,n}\mu_n^\bY(1,x_2,\dots,x_n)$, i.e. the probability that $\bY^v=1$ at some position $v$ on $\Sq(n)$ is, up to change of variables, equal to $\bar{\G^{\Sq(n)}}(x,y)$. A nice description of $\mu^{\bY}_{n}$ would thus be helpful. For this, two new ideas arise: \\
\tb The first one is to search for a representation of the same type as \eref{eq:prodd}, but for matrices $A_{x,y}$ (taking  the trace afterwards) instead of ``real numbers $A_{x,y}$''. This is not just a way to add some degrees of freedom: Proposition \ref{pro:rep} below says that all measures  invariant under rotation on $E_n$ have such a representation, and moreover have a representation of the form $\Trace(\prod_{i=0}^{n-1} Q^{x_i})$ for some matrices $Q^0$ and $Q^1$, \\
\tb The second idea comes from the product form \eref{main:equa1} of the transition kernel $\bT_n^\bZ$ for $\bZ \in \{\bX,\bY\}$. It suggests that a ``local equation'' linking the matrix $Q^x$ ``weighting the probability to observe $x$ somewhere on $\row(l)$'' and $Q^y$, $Q^{y'}$, the matrices weighting their two children in $\row(l+1)$ could suffice. But two neighbours in $\row(l)$ share a neighbour in $\row(l+1)$, and this must be taken into account. The second idea is to share the neighbour, its suffices to ``split'' the matrices, and search for a solution of the form $Q^x=V^xH^x$ for some matrices $V^x$ and $H^x$ having right splitting properties. This is what is done and proved to be possible, up to taking infinite matrices $Q^x$.

\subsection{Transfer matrix and measure representations}
\label{note:o}

We recall the definition and some properties of the Kronecker product (or tensor product) between matrices that we will use in the paper.  If $A$ 
%=\begin{bmatrix}
%A_{1,1}&\dots&A_{1,k} \\
%\vdots & \ddots&\vdots\\
%A_{m,1}&\dots&A_{m,k}
%A_{i,j}\end{bmatrix}_{1\leq i\leq m,j}$ 
is an $m\times k$ matrix and $B$ is a $p\times q$, then the Kronecker product $A\otimes{B}$ is the $mp\times kq$ matrix
 \[{A}\otimes{B} = \begin{bmatrix} A_{1,1} B & \cdots & A_{1,k}B \\ \vdots & \ddots & \vdots \\ A_{m,1} B & \cdots & A_{m,k} B \end{bmatrix}.\] 
The Kronecker product is associative. For matrices $A,B,C,D$, 
\ben 
\label{eq:t1}(A\otimes B) (C \otimes D)= (AC) \otimes (BD),\\
\label{eq:t2}\Trace(A\otimes B)=\Trace(A)\Trace(B)
\een
where it is assumed in \eref{eq:t1} that $A,B,C,D$ have sizes such that $AC$ and $BD$ are well defined, and in \eref{eq:t2} that the matrices are square matrices. The Kronecker product extends to infinite matrices $A\otimes B$, but it is ``interesting'' only for finite $B$. \par
For $n\geq 1$, let ${\cal M}_n$ be the set of probability measures on $E_n$, invariant under rotation: if $\mu\in{\cal M}_n$, then for any $\bx=(x_0,\dots,x_{n-1})\in E_n$, $\mu(\bx)=\mu(x_1,\dots,x_{n-1},x_0)$. 
\begin{pro}\label{pro:rep}
For any $\mu\in{\cal M}_n$, there exists two square matrices $Q^0$ and $Q^1$ of finite size, and four rectangular finite matrices  $V^0,V^1,H^0,H^1$  such that:\\
$(1)$ for any $\bx\in E_n$,
\[\mu(\bx)=\Trace\l(\prod_{i=0}^{n-1} Q^{x_i}\r),\]
$(2)$ $Q^1=V^1H^1,~~ Q^0=V^0H^0,~~ V^1H^0=\matz$ and $V^0H^1=\matz,$ where $\matz$ stands for the null matrix with the appropriate size.
\end{pro}

\noindent \Proof First, assume that there exists  some matrices $(Q^x, x\in\{0,1\})$ of size $m\times m$, such that  $\mu(\bx)=\Trace(\prod_{i=0}^{n-1} Q^{x_i})$.  Take  $\tilde V^1=Q^1\otimes \begin{bmatrix}1&0 \end{bmatrix}$, $\tilde H^1=\Id(m) \otimes \begin{bmatrix}1\\0 \end{bmatrix}$, $\tilde V^0=Q^0\otimes \begin{bmatrix}0&1 \end{bmatrix}$ and $\tilde H^0=\Id(m) \otimes \begin{bmatrix}0\\1 \end{bmatrix}$, where $\Id(m)$ is the identity matrix of size $m\times m$. Using $\Trace(\prod a_i\otimes b_i)=\Trace(\prod a_i) \Trace(\prod b_i)$, we see that for any $\bx\in E_n$,
\[\Trace\l(\prod_{i=0}^{n-1} Q^{x_i}\r)= \Trace\l(\prod_{i=0}^{n-1} \tilde{Q}^{x_i}\r),\] and clearly
$\tilde{Q}^1=\tilde{V}^1\tilde{H}^1,\tilde{Q}^0=\tilde{V}^0\tilde{H}^0, \tilde{V}^0\tilde{H}^1=\tilde{V}^1\tilde{H}^0=\matz.$

Hence, only (1) remains to be proved. For this, denote by ${\cal M}_n'$ the set of  probability distributions having the representation $\mu(\bx)=\Trace(\prod_{i=0}^{n-1} Q^{x_i})$ for some finite matrices $Q^0$ and $Q^1$. We show that ${\cal M}_n'={\cal M}_n$. For this, observe that ${\cal M}_n'$ is closed under finite mixture: if $\mu=\alpha\mu_a+(1-\alpha)\mu_b$ and  $\mu_j(\bx)= \Trace(\prod_{i=0}^{n-1} Q_j^{x_i})$ for $j\in\{a,b\}$, then taking
\[Q^x:=\begin{bmatrix}\alpha^{1/n}Q_a^x&0\\0&\beta^{1/n}Q_b^x  \end{bmatrix},~~~~x\in\{0,1\},\]
we obtain $\mu(\bx)= \Trace(\prod_{i=0}^{n-1} Q^{x_i})$. This shows that ${\cal M}_n'$ is closed under finite mixture. \par Since $\Trace(AB)=\Trace(BA)$, any probability distribution in ${\cal M}_n'$ is invariant under rotation. It remains to see that ${\cal M}_n'$ contains the uniform distribution on ``simple'' classes of rotation in $E_n$, since any measure in ${\cal M}_n$ is a mixture of these distributions. Take an element ${\bf \alpha}=(\alpha_0,\dots,\alpha_{n-1})$ in $E_n$. For any $i\in\cro{1,n}$, let ${\bf \alpha}(i)=(\alpha_i,\dots,\alpha_{n-1},\alpha_0,\dots,\alpha_{i-1})$ and  $R:=\{{\bf \alpha}(i),i\in\{1,\dots,n\}\}$, the rotation class of  ${\bf \alpha}$. We now show that the probability measure $\#R^{-1}\sum_{\beta\in R}\delta_{\bf \beta}$ belongs to ${\cal M}_n'$, which completes the proof. For this, 
denote by $z_i$ the number whose binary expansion is ${\bf \alpha}(i)$, and define $Q^1$ and $Q^0$ to be $2^n\times 2^n$ matrices where all entries are 0 except for
\[Q^{\alpha_i}_{1+z_i,1+z_{i+1}}=\#R^{-1/n},~~ \textrm{ for any }i.\]
It is then simple to check that $\mu(\bx)=\Trace(\prod_{i=1}^n Q^{x_i})$ coincides with $\#R^{-1}\sum_{\beta\in R}\delta_{\bf \beta}$.~$\Box$

\subsection{The fundamental Lemmas}
\label{sec:tosol }
The following two lemmas \ref{lem:finite} and \ref{pro:dec} -- that are among the main contributions of the present paper -- 
serve two different purposes. The first Lemma can be used to identify the solution  $\mu(\bx)=\Trace(\prod_i Q^{x_i})$ (and/or to prove that such a measure is a solution) of $\mu=\mu\bT_n$ when $\bT_n$ has a product form. It gives a sufficient condition on $(Q^0,Q^1)$ in the form of a ``finite'' system of quadratic equations. In the generic case, this system of equations does not depend on $n$, and then a right pair $(Q^0,Q^1)$ will provide a representation of $\mu$ for any $n$.\par
The second Lemma can be used when no such solution has been found. It permits us to  describe the measure $\nu \bT^\kappa$  on the $k$th line starting from $\nu$, for some particular $\nu$ using some matrices $\bQ^x_{(\kappa)}$. Since $\nu \bT^\kappa\sous{\to}{\kappa\to +\infty}\mu$ the solution of $\mu=\mu \bT$, 
this provides a way to approach $\mu$. The problem is that the matrices $\bQ_{(\kappa)}^x$ appear to grow with $\kappa$. The discussion of the convergence of the sequence $(\bQ_{(\kappa)}^x,\kappa)$ is addressed later in this paper: this gives the clues we mentioned in the abstract that $\G(x,y)$ should have an expression using some eigenvectors of some explicit infinite matrices.\par

Lemmas \ref{lem:finite} and \ref{pro:dec} are stated in the case of $\Sq(n)$ for processes with values in $\{0,1\}$. They will adapted to the triangular lattice and to processes with more than 2 values in Section \ref{sec:tri-lat}. \par

A matrix $\bT_n$ indexed by $E_n$ is called a transfer matrix. If for all $\ba\in E_n$, $\sum_{\bb} \bT_n(\ba,\bb)=1$, then $\bT_n$ is a probability transfer matrix. Moreover if, for any $\ba,\bb \in E_n$,
\beq\label{eq:prodform}
\bT_n(\ba,\bb)=\prod_{i=0}^{n-1} T_{a_i,a_{i\oplus 1},b_i}
\eq
for some $T_{a,b,c}$, $a,b,c,\in\{0,1\}$, we say that the matrix transfer has a product form.
\begin{lem}\label{lem:finite} Assume that $\bT_n$ is a probability transfer  matrix with a product form. Assume that there exists square matrices $(V^x, H^x,x\in\{0,1\})$ such that, $V^x H^{y}=\begin{bmatrix}0 \end{bmatrix}$ for $x\neq y$, and such that for any $x$
\beq \label{eq:finite}
V^x H^x=\sum_{y,y'\in \{0,1\}}  H^y V^{y'} T_{y,y',x}.
\eq
Finally let  $Q^x=V^{x}H^{x}$  for $x\in\{0,1\}$, and $\mu(\bx)=\Trace\l(\prod_{i=0}^{n-1} Q^{x_i}\r)$. We have $\mu=\mu\bT_n$.
\end{lem}
\begin{note}(1)
Equation \eref{eq:finite} is the system of quadratic equations mentioned in the introduction and in the Abstract. To find a solution to $\mu=\mu \bT_n$ it suffices to find $(V^x,H^x,x\in\{0,1\})$ satisfying \eref{eq:finite}.\\
(2) One may weaken the conditions in Lemma \ref{lem:finite}, replacing  \eref{eq:finite} by ``there exists an invertible matrix $P$'', such that 
\beq \label{eq:finite-2}
V^x H^x=P\l(\sum_{y,y'\in \{0,1\}}  H^y V^{y'} T_{y,y',x}\r)P^{-1}.
\eq
%Condition \eref{eq:finite} corresponds to the case where $P$ is the identity matrix.  
\end{note}
Before proving this Lemma we make two remarks. Firstly, we are only interested in non zero solutions! Under the hypothesis of this Lemma, the positivity of the measure $\mu$ is not guaranteed, nor the fact that it is non-zero or real. In the case where a non trivial solution exists, in the sense where $\mu(\bx)\neq0$ for some $\bx\in E_n$, then $\mu$ is a multiple of an eigenvector of $\bT_n$ with eigenvalue 1. Since $\bT_n$ is a Markov probability kernel, by the Perron-Frobenius theorem there is exactly one such eigenvector, which has non-negative entries. Note that if $(V^x,H^x,x\in \{0,1\})$ is a solution of \eref{eq:finite} then so is $(cV^x,H^x,x\in \{0,1\})$ for any $c\neq 0$. A $c$ (which may depend on $n$) can be taken such that $\mu$ is exactly the invariant distribution. Hence, a condition for non triviality (CNT) for a solution of \eref{eq:finite} is as follows:
\beq\label{eref:CNT1}
CNT_1:~~\textrm{ 1 is an eigenvalue of }\sum_{x\in{0,1}} V^xH^x.
\eq\\
Secondly, finding a $(V^x,H^x, x\in \{0,1\})$ solution of \eref{eq:finite} is difficult in general. The number of variables and of equations quadratically increase with the size of the matrices $(V^x,H^x,x\in\{0,1\})$.\medskip

\noindent\Proofof{ Lemma \ref{lem:finite}} Write $\Trace\l(\prod_{i=0}^{n-1} Q^{x_i}\r)=\Trace\l(\prod_{i=0}^{n-1} \sum_{\forall i, y_i,y'_{i}\in \{0,1\}}  H^{y_i} V^{y'_{i\oplus1}} T_{y_i,y'_{i\oplus1},x_i}\r)$.
Since $V^{y'_{i}}H^{y_{i}}=0$ if $y'_{i}\neq y_{i}$, we impose $y'_{i}=y_{i}$. It remains  to write, by some usual commutations 
\[\sum_{\by \in E_n}\Trace\l(\prod_{i=0}^{n-1}  H^{y_i} V^{y_{i\oplus1}}\r)\prod_{i=0}^{n-1} T_{y_i,y_{i\oplus1},x_i}=
\sum_{\by \in E_n}\mu(\by)\bT_n(\by,\bx)=\mu \bT_n.~\Box\]

We now state the second fundamental lemma.
\begin{lem}\label{pro:dec} Let $\bT_n$ be a transfer matrix on $E_n$. Consider two matrices $Q^x_{(0)}, x\in\{0,1\}$, of size $m\times m$ and  some  rectangular matrices $V^x_{(0)}, x\in\{0,1\}$ of size $m\times l$, and  $H_{(0)}^x, x\in\{0,1\}$ of size $l\times m$, such that
\[Q^x_{(0)}=V^x_{(0)}H^x_{(0)}, \textrm{ and, for  }x\neq y,~~V^x_{(0)}H^y_{(0)}=[0].\]
Consider the measure $\mu_{(0)}$ defined on $E_n$ by
$\mu_{(0)}(\bx)=\Trace\l(\prod_{i=1} Q_{(0)}^{x_i}\r).$ Assume that there exist 8 matrices $(h_{xy},(x,y)\in\{0,1\}^2)$ and  $(v_{xy},(x,y)\in\{0,1\}^2)$ (of size $i\times j$ and $j\times i$ so that the products $h_{xy}v_{x'y'}$ and $v_{xy}h_{x'y'}$ are defined) such that for any $\by$, $\bx$ in $E_n$,
\begin{equation}\label{eq:f1}
\Trace\l(\prod_{i=0}^{n-1} h_{y_i,x_{{i\oplus1}}}v_{y_{{i\oplus1}},x_{{i\oplus1}}}\r)= \bT_n(\by,\bx).
\end{equation}
Now, set for  $x \in \{0,1\}$ 
\beq\label{eq:reecriture}
V_{(1)}^x=\sum_{y\in\{0,1\}} H_{(0)}^y \otimes h_{y,x},~~H_{(1)}^x=\sum_{y\in\{0,1\}} V_{(0)}^y \otimes v_{y,x},~~Q_{(1)}^x=V_{(1)}^xH_{(1)}^x
\eq
and $\mu_{(1)}(\bx)=\Trace\l(\prod_{i=0}^{n-1} Q_{(1)}^{x_i}\r).$ Then under these conditions,\\
1)  $\mu_{(1)}=\mu_{(0)}\bT_n$. If moreover $\mu_{(0)}$ is a probability measure and $\bT_n$ is a probability transfer matrix, then  $\mu_{(1)}$ is also a probability measure. \\
2) moreover if for any $x,x'$,
\begin{equation}\label{eq:cond2}
h_{xy}v_{x'y'}=\begin{bmatrix}0 \end{bmatrix} \textrm{ for }y\neq y',
\end{equation} 
then $V^x_{(1)}H^y_{(1)}=\begin{bmatrix}0 \end{bmatrix}$ for $x \neq y$. \\
3) Assume that  the  probability transfer matrix $\bT_n$ has the product form \eref{eq:prodform} and that there exist $h_{x,y}$ with one line, and $v_{x,y}$  with one column such that
\begin{equation}\label{eq:tra}
h_{y,x}v_{y',x}=T_{y,y',x}, \textrm{ for any }y,y',x.
\end{equation}
Then condition \eref{eq:f1} is satisfied, and $Q_{(1)}^x$ and $Q_{(0)}^x$ have the same size. If moreover $Q_{(1)}^x=Q_{(0)}^x$ for $x\in\{0,1\}$, then  $\mu_{(1)}=\mu_{(0)}$.% [In fact, for any $T$ such that $\sum_cT_{a,b,c}=1$ for any $a,b\in\{0,1\}$ it is possible to find such mono-line and mono-column vectors $h_{a,b},v_{a,b}$]. 
\end{lem} 

\subsubsection*{Important note}
$\bullet$  The size of the cylinder $n$ plays no role at all in the hypothesis. Hence, the same matrices $H^x_{(j)}$ , $V^x_{(j)}$, $Q^x_{(j)}$ will (or will not) work for all $n$.   \par
$\bullet$   In the case where $Q_{(1)}^x=Q_{(0)}^x, x\in\{0,1\}$, we say that $Q_{(0)}$ is a fixed point for the matrix equation.  That is, for any $n$, $\mu_n$ defined on $E_n$ by
\[\mu_n(\bx)=\Trace\l(\prod_{i=0}^{n-1} Q_{(0)}^{x_i}\r),\] 
is a solution of $\mu_n=\mu_n\bT_n$.\medskip

Clearly, condition (2) is the condition needed to iterate the construction, that is to be able to give a matrix representation of the measure on the $\kappa$'th line, as stated in the next Corollary.
\begin{cor}\label{cor:big-construction} 
Let $\mu_{(0)}$, $Q_{(0)}^x$, $V_{(0)}^x$, and $H_{(0)}^x$ satisfy the hypothesis of Lemma \ref{pro:dec} (for example $V_{(0)}^1=H_{(0)}^1:=\begin{bmatrix}1 \end{bmatrix}$  and $V_{(0)}^0=H_{(0)}^0:=\begin{bmatrix}0 \end{bmatrix}$), and $h_{x,y}$ and $v_{x,y}$ satisfy \eref{eq:f1}, \eref{eq:cond2} and \eref{eq:tra}. Then, for any $\kappa\geq 1$, $x\in\{0,1\}$, let
\beq\label{eq:it-cons}
\left.\begin{array}{l}V_{(\kappa)}^x=\sum_{y\in\{0,1\}} H_{(\kappa-1)}^y \otimes h_{y,x},\\
H_{(\kappa)}^x=\sum_{y\in\{0,1\}} V_{(\kappa-1)}^y \otimes h_{y,x},
\end{array}\r.
\eq
and $Q^{x}_{(\kappa)}:=V_{(\kappa)}^xH_{(\kappa)}^x,~~ x\in\{0,1\}$, 
and 
\begin{equation}\label{eq:mu}
\mu_{(\kappa)}(\bx)=\Trace\l(\prod_{i=0}^{n-1} Q^{x_i}_{(\kappa)}\r)
\end{equation}
Then, for any $\kappa\geq 1$, $\mu_{(\kappa)}=\mu_{(\kappa-1)}\bT=\mu_{(0)}\bT^\kappa$. 
\end{cor}
\noindent\Proofof{Lemma \ref{pro:dec}} This proof relies on \eref{eq:t1} and \eref{eq:t2}; we also make heavy use of $\Trace(A_1\dots A_j)=\Trace(A_jA_1\dots A_{j-1})$.  For $(1)$, write
\be
\Trace\l(\prod_{i=1} Q_{(1)}^{x_i}\r)&=&\Trace\l(\prod_{i=1} \l(\sum_y H_{(0)}^y \otimes h_{y,x_i}\r)  \l(\sum_{y'} V_{(0)}^{y'} \otimes v_{y',x_i}\r)\r)\\
&=&\Trace\l(\prod_{i=1} \l(\sum_{y'} V_{(0)}^{y'} \otimes v_{y',x_i}\r) \l(\sum_y H_{(0)}^y \otimes h_{y,x_{i\oplus1}}\r) \r)\\
&=&\Trace\l(\prod_{i=1} \l(\sum_{y'} V_{(0)}^{y'} \otimes v_{y',x_i}\r) \l( H_{(0)}^{y'} \otimes h_{y',x_{i\oplus1}}\r) \r)\\
&=& \sum_{\by}\Trace\l(\prod_{i=1} \l(V_{(0)}^{y_i} \otimes v_{y_i,x_i}\r) \l(H_{(0)}^{y_i} \otimes h_{y_i,x_{i\oplus1}}\r) \r)\\
&=&\sum_{\by}\mu_{(0)} (\by)  \Trace\l(\prod_{i=0}^{n-1} v_{y_i,x_i}h_{y_i,x_{i\oplus1}}\r)\\
&=&\sum_{\by}\mu_{(0)} (\by)  \Trace\l(\prod_{i=0}^{n-1} h_{y_i,x_{i\oplus1}}v_{y_{i\oplus1},x_{i\oplus1}}\r)
\ee
 $(1)$ now follows from \eref{eq:f1}. \\
For (2), write for $x\neq y$, 
\be
V^x_{(1)}H^y_{(1)}&=&\l(\sum_{w} H_{(0)}^w \otimes h_{w,x}\r)\l(\sum_z V_{(0)}^z \otimes v_{z,y}\r)\\
&=&\sum_{w}\sum_z (H_{(0)}^w V_{(0)}^z) \otimes (h_{w,x}v_{z,y})
\ee
which is indeed $\matz$, since $h_{w,x}v_{z,y}=\matz$ for any $w,z$. $(3)$ follows immediately.~$\Box$

\subsection{Relation with ASEP/PASEP/TASEP}
The system \eref{eq:finite} is quite close to the systems associated with the ASEP/PASEP/TASEP problems which use the ``matrix ansatz'' proposed by  B. Derrida \& al. \cite{DEHP}. In this seminal paper, it is shown that the invariant distribution of the ASEP can be expressed as the formal result of a computation in some algebraic structure, where some operators $D$ and $E$ satisfy some quadratic relations of the form $DE=ED + E +D$ (or  more generally with some multiplicative parameters added) and some additional ``border conditions''of the type $D|V\!\!>=\beta|V\!\!>$,
 $<\!\!W|E=\alpha^{-1}<\!\!W|$. For these problems, matrix solutions  $D,E,V,W$  are explicitly found. These solutions have, depending on the values of the parameters and border conditions, a finite or an infinite size. In subsequent studies of these exclusion processes, finite/infinite matrices  appear  as solutions of quadratic equations of a type appearing in \cite{DEHP}. Each time, the questions of commutation of matrices, tensor products, and existence of limits with some growing matrices arise. We send the interested reader to \cite{BE} and \cite{HP}, and references therein. \par
These quadratic problems/algebras are also at the core of several  different problems of enumerative combinatorics, as observed by Viennot (this is discussed at length in several of his talks and courses, available on his web page, see e.g. \cite[section 7]{VG2}).

\section{Computation of fixed point solutions for gases of type 1 and 2}

As explained in the previous subsection, two main cases emerge. First, there may exist a fixed point for the matrix equation involving finite matrices according to Lemma \ref{pro:dec}(3). If no such fixed point solution exists (or even if it does), Corollary \ref{cor:big-construction} allows one to build bigger and bigger matrices to describe the distribution on the $k$th line. In some cases, taking the limit give rise to some infinite matrices. We first discuss, what happens when matrices of finite size are fixed points of the matrix equation. This is the case for the gas of type 1 on $\Sq(n)$.

\subsection{Gas of type 1}

This happens in the case $\bZ=\bX_{\Sq(n)}$, for which numerous different computations of the gas density exist. Let us add one way to find the solution.

\subsubsection{Application of Lemma \ref{lem:finite}}
\label{sec:al}  It suffices to search for finite matrices $(V^x,H^x,x\in\{0,1\})$ which solves \eref{eq:finite}.   Take 
\begin{equation}\label{eq:alpha}
V_{(0)}^0= \left[ \begin {array}{cc} 0&s_{{1,2}}\\ \noalign{\medskip}0&s_{{2,2}}
\end {array} \right], 
V_{(0)}^1:= \left[ \begin {array}{cc} s_{{1,1}}&0\\ \noalign{\medskip}s_{{2,1}}&0
\end {array} \right],  
H_{(0)}^0:=\left[ \begin {array}{cc} 0&0\\ \noalign{\medskip}t_{{1,2}}&t_{{2,2}}
\end {array} \right], 
H_{(0)}^1:= \left[ \begin {array}{cc} t_{{1,1}}&t_{{1,2}}\\ \noalign{\medskip}0&0
\end {array} \right].
\end{equation} 
They are chosen to trivially satisfy $V_{(0)}^xH_{(0)}^y=\matz$ if $x\neq y$. 
\begin{note}The  choice of the letters $V$ and $H$ comes from these vertical and horizontal structures.
\end{note} 
%The condition $H_{(0)}^x={}^tV_{(0)}^x$ can be arbitrary added d is taken to induce for free $\Trace(\prod_i Q_n^{x_i})=\Trace(\prod_i Q_n^{x_{n-i+1}})$, since we know that the measure $\mu$ solution of $\mu=\mu\bT^\bX$ owns this symmetry.  
The system of equations $V^xH^x=\sum_{y,y'} H^y V^{y'}\bT^\bX_{y,y',x}$ is equivalent to:
\[\Sys_1^\bX:=\left\{\begin{array}{l}
 s_{1,1}t_{2,1}=0,s_{2,1}t_{1,1}=0,t_{1,1}s_{1,1}=0,s_{2,2}t_{2,2}= ( 1-p )( s_{1,2}t_{1,2}+s_{2,2}t_{2,2}),\\
t_{2,1}s_{2,1}=p ( s_{1,2}t_{1,2}+s_{2,2}t_{2,2} ) ,s_{1,2}t_{1,2}=t_{1,1}s_{1,1}+t_{2,1}s_{2,1},\\
s_{1,2}t_{2,2}=t_{1,1}s_{1,2}+t_{2,1}s_{2,2},s_{2,2}t_{1,2}=t_{1,2}s_{1,1}+t_{2,2}s_{2,1} 
\end{array}\right.\]
This system has non trivial solutions. For example
$s_{1,1}=0,s_{1,2}=p,s_{2,1}= 1/( 1-p )
,s_{2,2}=1,t_{1,1}=0,t_{1,2}=1,t_{2,1}= ( 1-p ) p
,t_{2,2}=1-p$,
 in which case
\[Q^0=\left[ \begin {array}{cc}p&p(1-p)
\\ 1&1-p\end {array} \right] ,Q^1= \left[ \begin {array}
{cc} 0&0\\ \noalign{\medskip}0&p\end {array} \right].\]
From here the density  $\Trace((Q^0+Q^1)^{n-1}Q^1)/\Trace((Q^0+Q^1)^n)$ of the gas of type 1 on the cylinder can be computed, and from that  $\G^{\Sq(n)}(x,1)$ thanks to Proposition \ref{pro:fo}. Taking the limit gives $\G^{\Sq}(x,1)$. This representation was known and is present under different forms in \cite{NDV,DH1,BM1}. What is remarkable is that it follows from a simple computation.

\subsubsection{Application of Lemma \ref{pro:dec}}

For the same result, one may use Lemma \ref{pro:dec} since a solution to both \eref{eq:tra} and  \eref{eq:cond2} exists.
The idea is to search for solutions under the following form, for which \eref{eq:cond2} automatically holds~:
\ben\label{eq:les-h}
\begin{array}{ccc}
h_{0,0}&=\begin{bmatrix} 0,a_1,0,a_2\end{bmatrix},~~
h_{1,0}&=\begin{bmatrix} 0,b_1,0,b_2\end{bmatrix}\\ 
h_{0,1}&=\begin{bmatrix} c_1,0,c_2,0\end{bmatrix},~~
h_{1,1}&=\begin{bmatrix} d_1,0,d_2,0\end{bmatrix}\\ 
\end{array}
\een
and $v_{x,y}={}^t h_{x,y}$ for all $x,y$,  this last condition not being needed. % but it implies the symmetry of a fix point solution $Q_{(1)}^\star=Q_{(0)}^\star$ symmetric.  
The system \eref{eq:tra} can be rewritten as~:
\beq\label{eq:sys1}
\Sys_2^\bX:=\left\{
\begin{array}{l}
d_1^2+d_2^2=0,~c_1d_1+c_2d_2=0,~b_1^2+b_2^2-1=0,\\
c_1^2+c_2^2-p=0,~b_1a_1+b_2a_2-1=0,~a_1^2+a_2^2-1+p=0.
\end{array}
\right.
\eq
Clearly, $d_1=0, d_2=0$, and the rest of the system becomes
\[\Sys^\bX:=\{b_1^2+b_2^2-1=0,~{c_{1}}^{2}+{c_{2}}^{2}-p=0,~{a_{1}}^{2}+{a_{2}}^{2}-1+p=0,~b_{1}a_{1}+b_{2}a_{2}-1=0\}.\]
Solutions to this system exist. %Indeed, consider 3 points  $A=(a_1,a_2)$, $B=(b_1,b_2)$, $C=(c_1,c_2)$ in $`R^2$; $\Sys^\bX$ says that $A,B,C$ lie on some circle with center 0, and radius $\sqrt{1-p},1,\sqrt{p}$, and that $A.B=1$; for   $p\in(0,1)$, there are no solutions: the existing solutions for $(a_1,a_2,b_1,b_2,c_1,c_2)$ are (at least partially) complex. 
We then take $V_{(0)}^x$, $H_{(0)}^x$ as in \eref{eq:alpha} for example, and take $V_{(1)}^x$,$H_{(1)}^x$ as defined in \eref{eq:it-cons}.  Again, $Q^x_{(1)}=Q^x_{(0)}$ is possible since $h_{x,y}$ and $v_{x,y}$ are respectively single line and single column (according to Lemma \ref{pro:dec}). Of course, 
since $h_{y,x}v_{y',x}=T^\bX(y,y',x)$, we are back to the previous problem treated in Section \ref{sec:al}. 
\begin{note}
The appearance of non real numbers in the considerations ($\Sys^\bX$ has no real solution) does not harm the reasoning at all  since  $\bT^\bY$ is still a probability transfer matrix. The possibility to use non real matrices $V^x,H^x$ as a solution of the equations of interest enriches the space of solutions.
\end{note}

\subsection{No finite solution for the gas of type 2}

In the case of the gas of type 2, we were not able to find finite (non-trivial) matrices $(V^x,H^x,x\in\{0,1\})$ which solve \eref{eq:finite}.  If they exist, they must have size $>5$ (this can be seen via the computation of a Gröbner basis). Also, numerical computations -- which in principle does not guarantee any result -- indicate that no non-trivial matrices $(V^x,H^x,x\in\{0,1\})$ with complex coefficients having size $\leq 8$ are solution. Hence we were unable to use Lemma \ref{lem:finite} to go on.
 \begin{Ques}(a) There it exist finite square matrices $V^0,H^0,V^1,H^1$ which solve \eref{eq:finite} for $T=T^\bY$? (and such that  $\sum_x V^xH^x$ has eigenvalue 1?)\\
(b) If not, is it possible to find a solution to \eref{eq:finite-2} for $T=T^\bY$?
 \end{Ques}

\section{Toward an infinite size solution}

This section prospective explains how some infinite matrices may arise. Even if no complete solution is provided, we hope that this new point of view will allow some reader to tackle the problem of computing $\G(x,y)$.\par

Recall  that Lemma \ref{pro:dec} and Corollary \ref{cor:big-construction} say that if there exist single line and single column matrices $h_{x,y}$ and $v_{x,y}$ solving $h_{x,y}v_{x',y'}=T^\bZ(x,x',y) 1_{y=y'}$, then on a cylinder of the square lattice the distribution on the $\kappa+1$'th line $\mu_{(\kappa)}:=\mu_{(\kappa-1)}\bT^\bZ$ has a representation of the form \eref{eq:mu}, provided that the distribution on the first line $\mu_{(0)}$ has this same form with $V_{(0)}^xH_{(0)}^y=0$ for $x\neq y$. Such matrices $h_{x,y}$ and $v_{x',y'}$ exist for $\bZ\in\{\bX,\bY\}$:
\begin{theo}\label{theo:mach} For any $\bZ\in\{\bX,\bY\}$,\\
(1) there exist single line matrices $(h_{x,y},x,y\in\{0,1\})$ and single column matrices $(v_{x,y},x,y\in\{0,1\})$ which solves \eref{eq:cond2} and \eref{eq:tra}, namely
\begin{equation}\label{eq:master-equation}
h_{x,y}v_{x',y'}=T^\bZ(x,x',y) 1_{y=y'} \textrm{ for any } (x,y,x',y')\in\{0,1\}^4.
\end{equation}
(2) For any $\kappa$, the entries of the matrices $V_{(\kappa)}^x,H_{(\kappa)}^x,Q_{(\kappa)}^x$ (as introduced in Corollary \ref{cor:big-construction}) can be computed. \\
(3)  The solutions of \eref{eq:master-equation} and the initial matrices  $V_{(0)}^x,H_{(0)}^x$ can be chosen in such a way that for any $x$, $V_{(2\kappa)}^x, H_{(2\kappa)}^x, Q_{(2\kappa)}^x$ converges simply when $\kappa\to +\infty$ (that is, each fixed entry converges).
\end{theo}

We have already discussed the existence of solutions to $\Sys_2^\bX$, this implying the existence of solutions to \eref{eq:master-equation} in the case $\bZ=\bX$. To prove Theorem \ref{theo:mach}(1), let us write the following system of equations for $h_{x,y}$ and $v_{x,y}$ defined as in \eref{eq:les-h}, in the case $\bZ=\bY$. In this case  \eref{eq:master-equation} is equivalent to
 \begin{equation}\label{eq:sys2}
 \Sys_2^\bY:=\left\{
 \begin{array}{l}
 {c_1}^2+{c_2}^2-q+pq=0,~{c_1}^2+q{d_1}^2-q{c_1}^2=1, {d_1}^2+{d_2}^2-p+pq-q=0 \\
 d_1c_1+d_2c_2-q+pq=0, ~{a_1}^2+{a_2}^2-1+q-pq=0,\\ 
~b_1a_1+b_2a_2-1+q-pq=0,~{b_1}^2+{b_2}^2-1+q+p-pq=0.
 \end{array}\right.
 \end{equation}
It is not difficult to check that this system has a solution (using Maple or Mathematica, for example, or the computation of a Gröbner basis).

The proofs of Theorem \ref{theo:mach} (2) and (3) are more delicate. Their respective proofs are the object of Sections \ref{v:const} and \ref{seq:conv} below.\medskip

\noindent\bf Notation. \rm For any pair of matrices $(A^0,A^1)$, $A^\star:=A^0+A^1$. Similarly, for any doubly indexed quantity $a_{x,y}$, $a_{x,\star}=a_{x,0}+a_{x,1}$.

\subsection{Computations of the entries of $Q_{(k)}^x, V_{(k)}^x, H_{(k)}^x$} 
\label{v:const}

Let $\bZ\in\{\bX,\bY\}$ be fixed, and let $h_{x,y}$ and $v_{x,y}$ as defined in \eref{eq:les-h} be solution of $\Sys_2^\bZ$ (the sizes of $h$ and $v$ are $1\times m$ and $m\times 1$  respectively  for some $m\geq 1$).
Let us compute the entries of $\l((Q^x_{(\kappa)},V^x_{(\kappa)}, H^x_{(\kappa)}), \kappa \geq 1\r)$ starting from  some matrices  $V^1_{(0)}$, $H^1_{(0)}$, $V^0_{(0)}$ and  $H^0_{(0)}$ such that $V_{(0)}^xH_{(0)}^{x'}=\matz$ for $x\neq x'$ of size $m_0\times m_0$ (for example, $V^1_{(0)}=H^1_{(0)}=[1], V^0_{(0)}=H^0_{(0)}=[0]$). \par
For $\kappa\geq 2$, by associativity of the Kronecker product,
\begin{equation}\label{eq:V}
V^x_{(\kappa)}=\sum_y \sum_z (V_{(\kappa-2)}^z\otimes v_{z,y}) \otimes h_{y,x}  = \sum_z V_{(\kappa-2)}^z\otimes  \bv_{z,x}
\end{equation}
with
 \[\bv_{z,x}= \sum_y v_{z,y}\otimes h_{y,x},\]
which entails, iteratively that \begin{equation}\label{eq:rew-W}
V_{(2\kappa)}^{z_\kappa}= \sum_{z_0,z_1,\dots,z_{\kappa-1}} V_{(0)}^{z_0} \otimes \bv_{z_0,z_1} \otimes \bv_{z_1,z_2}\dots\otimes \bv_{z_{\kappa-1},z_{\kappa}}.
\end{equation}
(The matrices $\bv_{z,y}$ have size $m\times m$.)  A similar formula exists for $H_{(2\kappa)}^x$, obtained by replacing  $h$ by $v$ and vice versa in the previous considerations, leading to the definition of $\bh_{z,x}= \sum_y h_{z,y}\otimes v_{y,x}$. \par
In the same manner, $Q_{(\kappa)}^x=V_{(\kappa)}^x H_{(\kappa)}^x$ can be computed: write
 \be
Q_{(\kappa)}^x&=&\l(\sum_z V_{(\kappa-2)}^z\otimes  \bv_{z,x}\r)\l(\sum_{z'} H_{(\kappa-2)}^{z'}\otimes  \bh_{z',x}\r).
\ee
Using the structure of $V_{(\kappa-2)}^z$ and $H_{(\kappa-2)}^{z'}$, we get
\be
Q^{x}_{(\kappa)} &=& \sum_{w} Q^w_{(\kappa-2)} \otimes \bq_{w,x}
\ee
where
\be
\bq_{w,x}&=&\sum_y \sum_z\l(v_{w,y}h_{w,z}\r) \otimes \l(h_{y,x}v_{z,x}\r)
=\sum_y \sum_z  T^\bZ_{y,z,x}\l(v_{w,y}h_{w,z}\r)
\ee
(we used here that $v_{a,b},h_{a,b}$ is a solution of $\Sys_2^\bZ$). Again, using the same methods above
\begin{equation}\label{eq:rew-Q}
Q_{(2\kappa)}^{z_\kappa}= \sum_{z_0,z_1,\dots,z_{\kappa-1}} Q_{(0)}^{z_0} \otimes \bq_{z_0,z_1} \otimes \bq_{z_1,z_2}\otimes\dots\otimes \bq_{z_{\kappa-1},z_{\kappa}}.
\end{equation}
With these formulas, the entries of  $V_{(2\kappa)}^x, H_{(2\kappa)}^x$ as well as that of $Q_{(2\kappa)}^x$ can be computed (and with a simple adaptation those of odd indices also).
For this, recall that if
\[A=B\otimes C\] 
where $C$ is a $c\times c$ matrix, $B$ a $b\times b$ matrix then, for any $i,j \in \{0,\dots, bc-1\}$,
\beq\label{eq:ABC-entries}
A[{i,j}]=B[{i \Div m, j \Div m}]\times C[{i \Mod m, j \Mod m}]
\eq
with the convention that for any matrix $M$, $M[i,j]:=M_{i+1,j+1}$, and as usual $x \Div m$ and $x \Mod m$ denote the quotient and the remainder in the division of $x$ by $m$. 

Assume now that some matrices $\bW_{(2\kappa)}^x,x\in\{0,1\}$, $\kappa\geq0$, satisfy
\beq\label{eq:WW}
\bW_{(2\kappa)}^{z_\kappa}= \sum_{z_0,z_1,\dots,z_{\kappa-1}} \bW_{(0)}^{z_0} \otimes \bw_{z_0,z_1} \otimes \bw_{z_1,z_2}\otimes\dots\otimes \bw_{z_{\kappa-1},z_{\kappa}}
\eq
with $\bw_{x,y}$ having size $m\times m$, and $\bW_{(0)}^x$ having size $m_0\times m_0$. Therefore $\bW^x_{(2\kappa)}$ has size $m_0 \times m^\kappa$; any $i,j$ in $\{0,1,\dots,m_0 \times m^\kappa\}$ can be written under the following form :
\[i=a_{\kappa+1}(i) m^\kappa+ \sum_{l=1}^\kappa m^{l-1}a_l(i),~~j=a_{\kappa+1}(j) m^\kappa + \sum_{l=1}^\kappa m^{l-1}a_l(j),\]
where $0\leq a_{\kappa+1}(i),a_{\kappa+1}(j) <m_0$, and $0\leq a_l(i),a_l(j) <m$ for $l \in \cro{1,\kappa}$ (apart from $a_{\kappa+1}$ which may play a special role if $m_0\neq m$, the $a_{l}(i)$'s are the digits of $i$ in base $m$). Therefore, from \eref{eq:ABC-entries},
\beq \bW_{(2\kappa)}^{z_\kappa}[i,j]=\sum_{z_0,\dots,z_{\kappa-1}} \bW_{(0)}^{z_0}[a_{\kappa+1}(i),a_{\kappa+1}(j)] \times \bw_{z_1,z_2}[a_{\kappa}(i),a_{\kappa}(j)]\times\dots\times \bw_{z_{\kappa-1},z_{\kappa}}[a_1(i),a_1(j)].
\eq
There is also a way to represent this with matrices, very similar to the representation of the distribution of a Markov chain; for any $(a,b)$ in $\{0,\dots,m-1\}^2$, let $M_\bw(a,b)$ be the $2\times 2$ matrix defined by:
\[M_\bw(a,b):=
\begin{bmatrix}
\bw_{0,0}[a,b]&\bw_{0,1}[a,b] \\
\bw_{1,0}[a,b]&\bw_{1,1}[a,b]
\end{bmatrix},\]
and for $a,b \in \{0,\dots,m_0-1\}$,
\[\rho_{\bW}[a,b]:=
\begin{bmatrix}
\bW_{(0)}^0[a,b]&\bW_{(0)}^1[a,b]
\end{bmatrix}.\]
We have 
\beq \label{eq:V-kappa}
\bW_{(2\kappa)}^z[i,j]=\rho_\bW[a_{\kappa+1}(i),a_{\kappa+1}(j)]\, M_{\bW}[a_{\kappa}(i),a_{\kappa}(j)]\dots M_{\bW}[a_{1}(i),a_{1}(j)]
\begin{bmatrix}
1_{z=0} \\
1_{z=1}
\end{bmatrix}.
\eq
This fact, together with \eref{eq:rew-Q} and \eref{eq:rew-W}, proves Theorem \ref{theo:mach}(2).

\subsection{Convergence of the entries of $V_{(2\kappa)}^x,H_{(2\kappa)}^x$, and $Q_{(2\kappa)}^x$}
\label{seq:conv}
We continue from the previous section. 
Notice that for a fixed $(i,j)$, $a_{l}(i)$ and all $a_l(j)$ are zero for large $l$. We then immediately have:
\begin{lem}\label{lem:conv1}  $\bW_{(2\kappa)}^z[i,j]$ converges when $\kappa\to +\infty$ in $\mathbb{C}$ for any $i,j$ in $\mathbb{N}$ 
if and only if 
\[\rho_\bW[0,0] M_\bw[0,0]^l\] 
converges when $l$ goes to $+\infty$; a sufficient condition is the convergence of $M_\bw[0,0]^l$.
\end{lem}
Again, the convergence stated in Lemma \ref{lem:conv1} will be only interesting if the limit is not zero. 
We examine the simple convergence of $V_{(2\kappa)}^x, H_{(2\kappa)}^x$ and $Q_{(2\kappa)}^x$ when $\kappa\to+\infty$, to some infinite matrices 
$(V_{\infty}^x, H_{\infty}^x,Q_{\infty}^x)$, meaning that, for $x\in \{0,1\}$ and any $i,j\geq 0$,
\beq
V_{(2\kappa)}^x[i,j]\to V_\infty^x[i,j],~~H_{(2\kappa)}^x[i,j]\to H_\infty^x[i,j]~~\textrm{ and }~~Q_{(2\kappa)}^x[i,j]\to Q_{\infty}^x[i,j]
\eq
starting with some suitable matrices $((V_{(0)}^x,H_{(0)}^x,Q_{(0)}^x),x\in\{0,1\})$. If such a convergence holds, 
the limiting infinite matrices $Q_{\infty}^1$ and $Q_{(\infty)}^0$ are moreover a solution of the following rewriting 
like system:
\beq \label{eq:rew1}
\textrm{ for }x\in\{0,1\},~~ 
Q^x_{\infty} = \sum_{w=0}^1 Q^w_{\infty} \otimes \bq_{w,x}.
\eq
Rewriting rules such as \eref{eq:rew1} rely entirely 
on the corners of the matrices $(Q^x_{\infty}[0,0],x\in\{0,1\})$. Writing $\rho_M[0,0]:=[M^0[0,0],M^1[0,0]]$, formula \eref{eq:rew1}
allows us to compute the only possible corners:
\be
\rho_{Q_\infty}&=&\rho_{Q_\infty}\bM_\bq[0,0]\\
%Q_{(2\kappa)}^1[0,0]&=&\rho_{Q}[0,0]\bM_\bq[0,0]\begin{bmatrix}0 \\ 1\end{bmatrix}= ( q+p-pq ) d_1^2\\
\rho_{V_{\infty}}&=&\rho_{V_{\infty}}\bM_\bv[0,0]\\
\rho_{H_{\infty}}&=&\rho_{H_{\infty}}\bM_\bh[0,0].
%V_{(2\kappa)}^1[0,0]&=&\rho_{V}[0,0]\bM_\bv[0,0]\begin{bmatrix}0 \\ 1\end{bmatrix}=1
\ee
%
%with the same choice for $\rho_Q=\rho_V=\begin{bmatrix}0,1 \end{bmatrix}$ as before. This fixes the only suitable choice for
% $(Q^x_{\infty}[0,0],x\in\{0,1\})$, and then for $(Q^x_{(2\kappa)}[0,0],x\in\{0,1\})$. 
% \begin{note}
% Instead of the convergence of $Q^x_{(2\kappa)}$ one may also investigate the convergence of $Q^x_{(2\kappa)}/\lambda^\kappa$ for some $\lambda$, since
% $\Trace(\prod_{i=1}^n Q^{x_i}_{(2\kappa)})/\Trace(\prod_{i=1}^n (Q^{0}_{(2\kappa)}+Q^{1}_{(2\kappa)})^n)$ is invariant by dividing the numerator and the denominator by $\lambda^{n\kappa}$. 
% \end{note}
We examine separately the two cases $\bZ=\bX$ and $\bZ=\bY$ in the next subsections.

\subsubsection{Case of the gas of type 1.}

We work with $h$ defined in \eref{eq:les-h}. We then find
\[\begin{array}{ll}
M_\bv(0,0)=\begin{bmatrix}
0&c_1d_1 \\
0&d_1^2
\end{bmatrix}, &
M_\bq(0,0)=\begin{bmatrix}
c_1^2&0 \\
d_1^2&0
\end{bmatrix}\end{array}.\]
Therefore, the convergence of the sequence $(V_{(2\kappa)}^z[i,j],\kappa >0)$ to a non zero limit is equivalent to $d_1=1$, in which case, for any $l\geq 1$, 
\begin{equation}\label{eq:Mcvv}M_\bv(0,0)^l= M_\bv(0,0)\end{equation}
and the convergence of $M_\bq(0,0)^l$  arises if $c_1^2=1$, in which case, for any $l\geq 1$,
\begin{equation}\label{eq:Mcvq}M_\bq(0,0)^l=M_\bq(0,0).\end{equation}
We may wonder if $\Sys_2^\bX$ still has some solutions if we add these conditions. \\
\tb The answer is yes for the condition $c_1=1$ (for example, $a_1^2=-p,a_2=1, b_1=0,b_2=1,c_1=1,c_2^2=p-1,d_1=d_2=0$). Then it is possible to have simple convergence for $Q_{(2\kappa)}^x$.\\
\tb The answer is no for the condition $d_1=1$. 
In order to find $h_{x,y},v_{x,y}$ satisfying \eref{eq:master-equation} (and $\bZ=\bX$) and such that $V_{(2\kappa)}$ simply converges, it suffices to increase the size of the matrices $h$ and $v$ defined in \eref{eq:les-h}. Take instead
\ben\label{eq:les-h2}
\begin{array}{ccc}
h_{0,0}&=\begin{bmatrix} 0,a_1,0,a_2,0,a_3\end{bmatrix},~~
h_{1,0}&=\begin{bmatrix} 0,b_1,0,b_2,0,b_3\end{bmatrix}\\ 
h_{0,1}&=\begin{bmatrix} c_1,0,c_2,0,c_3,0\end{bmatrix},~~
h_{1,1}&=\begin{bmatrix} d_1,0,d_2,0,d_3,0\end{bmatrix}\\ 
\end{array}\een
and again $v_{x,y}={}^t h_{x,y}$. The values of $M_\bq(0,0)$ and  $M_\bv(0,0)$ are unchanged, but this time there are some solutions for \eref{eq:f1} and \eref{eq:cond2} and where $d_1=1$. Again, this is not difficult to check with a program like Maple or Mathematica. Note also that if we want to solve  $\Sys_2^\bX$ together  with the two equations $c_1=1$ and $d_1=1$, there exists solutions for $h$ and $v$ having size 6. For example:
\[a_1^2+p=0,a_2=0,a_3=1,b_1=0,b_2=0,b_3=1,c_1=1,c_2^2=p,c_3=i,d_1=1,d_2=0,d_3=i.\]
It remains to specify $\rho_\bq$ and $\rho_\bv$, namely the starting condition of the construction. We may take $m_0=1$ (that is, starting with $1\times 1$ matrices). Taking $Q^{(0)}=V^{(0)}=H^{(0)}=[0]$ and $Q^{(1)}=V^{(1)}=H^{(1)}=[1]$  leads to 
\beq \label{eq:start}
\rho_V=\rho_H=\rho_Q:=  \begin{bmatrix}0,1  \end{bmatrix}.%,~~ \rho_Q=\begin{bmatrix}0,1 \end{bmatrix}.
\eq 
With this convention,  by \eref{eq:Mcvv} and \eref{eq:Mcvq}, the upper-left corner of $Q_{(2\kappa)}^x$ coincides with $Q_{(2\kappa-2)}^x$, and the same thing hold for $V$ and $H$ as well. We then have $\rho_{Q_\infty}=\rho_{H_{\infty}}=\rho_{V_{\infty}}=\rho_{Q}$.

\subsubsection{Case of the gas of type 2.} 
First, for $h$ and $v$  defined in \eref{eq:les-h},
\[M_\bq(0,0)= 
\begin{bmatrix}  ( -1+p) ( -1+q ) c_1^2& ( q+p-pq ) c_1^2\\
 ( -1+p )  ( -1+q ) {d_1}^2& ( q+p-pq ) d_1^2 
\end{bmatrix},~~ 
M_\bv(0,0)=\begin{bmatrix}  0&d_1c_1\\ 0&{d_1}^2\end{bmatrix} .\]
The convergence of $M_\bq(0,0)^l$ to a non zero limit happens if $( p-pq+q ) {d_1}^2+ ( 1-q-p+pq ) c_1^2=1$, in which case 
$M_\bq(0,0)^l= M_\bq(0,0)$ for any $l\geq1$, and the convergence of  $M_\bv(0,0)^l$ to a non zero limit arises if $d_1=1$, in which case 
$M_\bv(0,0)^l= \begin{bmatrix}  0&c_1\\ 0&1
\end{bmatrix} $ for $l\geq 1$. In this case, solutions  $h_{x,y}$  of size  $1\times 4$ exist : there exists solutions to $\Sys_2^\bY$ with the additional condition $( p-pq+q ) {d_1}^2+ ( 1-q-p+pq ) c_1^2=1$ or  $d_1=1$.  
Again, if if we want to solve  $\Sys_2^\bY$ with both conditions together there exist solutions for $h$ and $v$ having size 6, for example:
\be
{b_3}^2&=& (1-q)(1-p), {d_2}^2=-(1-q)(1-p),\\
{a_1}^2 & =& -\frac{p(1-q+pq)}{(1-p)(1-q)} ,a_2=0,a_3=\frac {1-q+pq}{b_3},b_1=0,b_2=0,c_1=1,\\
c_2&=&-\frac {1-q+pq}{d_2}, c_3^2= \frac{p(1-q+pq)}{(1-p)(1-q)} ,d_1=1,d_3=0.
\ee
This suffices to imply the simple convergence of $Q_{(2\kappa)}^x, V_{(2\kappa)}^x, H_{(2\kappa)}^x$. Here $\rho_{Q_{\infty}}=[(-1+p)(-1+q),p-pq+q]$, $\rho_{V_{\infty}}=[0,1]$ and $\rho_{H_{\infty}}=[0,1]$.\\
This ends the proof of Theorem \ref{theo:mach}(3).

\subsection{Trace of the limit and limit of the trace}

As said above, Theorem \ref{theo:mach} gives a representation of $\mu_{(\kappa)}^\bZ$ starting from some simple $\mu_{(0)}^\bZ$ (this could be useful to make advances on enumeration issues concerning DA with height $\kappa$).
The important and natural question is the following one : do we have, for any $\bx=(x_1,\dots,x_k)$,
\beq
\mu_{(2\kappa)}(\bx)=\Trace\l(\prod_{i=0}^{n-1}Q_{(2\kappa)}^{x_i}\r)\xrightarrow[\kappa\to\infty]{}\mu_{\infty}(\bx)=\Trace\l(\prod_{i=0}^{n-1}Q_{\infty}^{x_i}\r) ? 
\eq
Since the simple convergence of $Q_{(2\kappa)}^x$ to $Q_{\infty}^x$ does not imply the simple convergence of $Q_{(2\kappa)}^xQ_{(2\kappa)}^y$ to $Q_{\infty}^x Q_{\infty}^y$,  the answer to this question is certainly not an immediate issue. For the TASEP, the choices of matrices $D,E$ satisfying the different matrix ansatz, leads or not to the convergence of the product (see discussions in \cite{DEHP,BE}).

Moreover, simple convergence of a sequence of matrices $A_n$ to some matrice $A_{\infty}$ does imply the convergence of the trace, since the trace involves an infinite number of entries. Nevertheless, $\mu_{(2\kappa)}$ converges when $\kappa$ goes to +$\infty$ by Lemma \ref{lem:MC}. Let $\mu^\infty(\bx)=\lim_\kappa \mu_{(2\kappa)}(\bx)$ the limit of the measure. The question is: do we have $\mu^\infty=\mu_\infty ?$ Since $\mu^\infty$ is the only non trivial probability measure fixed point of $\mu^\infty= \mu^\infty\bT^\bY$, it suffices to show that $\mu_\infty$ satisfies the same property, which would imply $\mu_\infty=\lambda\mu^\infty$, for some $\lambda$ (which must be shown to be $\neq 0$). 
In fact, by construction $q_{y,x}$ is associated with two-row transitions, since
$\Trace(\prod_{i=1}^n q_{y_i,x_i})=\sum_{\bz \in E_n} \prod T^\bY_{y_i,y_{i\oplus 1},z_i}T^\bY_{z_i,z_{i\oplus1},x_i}=(\bT^\bY)^2(\bx,\by)$. Clearly $(\bT^\bY)^2$ is also a  probability transfer matrix, corresponding to an aperiodic irreducible Markov chain on a finite state space. Subsequently, 
by uniqueness, it is easy to check that $\mu_\infty$ is the unique solution to  $\mu = \mu(\bT^\bY)^2$.

Using that $Q_{\infty}$ is solution of the rewriting system \eref{eq:rew1}, if one ignores convergence and commutation issues, then 
\ben\label{eq:mult-fix-point}
\Trace\l(\prod_{i=0}^{n-1}Q_{\infty}^{x_i}\r)&=&\sum_{\by \in E_n} \Trace\l(\prod_{i=0}^{n-1}Q_{\infty}^{y_i}\r)\Trace\l(\prod_{i=0}^{n-1} q_{y_i,x_i}\r)\\
&=& \sum_{\by \in E_n} \Trace\l(\prod_{i=0}^{n-1}Q_{\infty}^{y_i}\r)(\bT^\bY)^2(\by,\bx),
\een
and one sees that we just have to justify convergence of $\prod_{i=0}^{n-1}Q_{\infty}^{x_i}$ and the validity of the rearrangements in the infinite sum.

We were unable to prove the validity of this. Besides, the entries of $Q^y_\infty$ seems no to converge to 0 (due to the choice of the value of $M_\bq(0,0)$, needed to have the convergence of $Q_{(\kappa)}^x$ to $Q_\infty^x$); also, seen as series in $p,q$, the degrees of the entries $Q^\star_{i,j}$ do not go to $+\infty$ with $i,j\to \infty$. 
\begin{Ques}
Is it possible to prove that $\Trace\l(\prod_{i=0}^{n-1}Q_{\infty}^{x_i}\r)$ is well defined (for some notion of convergence) and solution of \eref{eq:mult-fix-point}?
\end{Ques}

We here review some properties of the matrices $V_{(\kappa)}^x,H_{(\kappa)}^x$ and of $\mu_{(2\kappa)}$. These properties lead to some questions about the structure of $Q_\infty^x$, and its eigenvectors (if any). First, we have
\be
\mu_{(\kappa)}(\bx)&=&\Trace\l(\prod_{i=0}^{n-1} V^{x_i}_{(\kappa)}H^{x_i}_{(\kappa)}\r)
=\Trace\l(\prod_{i=0}^{n-1} H^{x_i}_{(\kappa)}V^{x_{i\oplus1}}_{(\kappa)}\r).
\ee
This representation is very close to the standard representation of Markov chain, where here $ H^{a}_{(\kappa)}V^{b}_{(\kappa)}$ plays the role of a probability transition and
\[\bP_{(\kappa)}:=\sum_{a,b}  H^{a}_{(\kappa)}V^{b}_{(\kappa)}\]
plays the role of the transition matrix.
\begin{pro} For any $\kappa$, $Q_{(\kappa)}^\star$ has eigenvalues 1 and 0. The eigenvalue 1 has multiplicity 1.
\end{pro}

\Proof We first claim that $(Q_{(2\kappa)}^\star)^{2\kappa}(1-Q_{(2\kappa)}^\star)=(Q_{(2\kappa+1)}^\star)^{2\kappa}(1-Q_{(2\kappa+1)}^\star)=\matz$. The claim implies that the minimal polynomial of $Q_{(2\kappa+`e)}^\star$ (with $`e\in\{0,1\}$) divides $x^{2\kappa}(1-x)$, which implies that the eigenvalues are 0 and 1. Since $\Trace(Q^\star_{(\kappa)})=\Trace(Q^\star_{(\kappa-2)})$, and since $\Trace(Q_{(0)}^\star)=\Trace(Q_{(1)}^\star)=1$, the eigenvalue 1 has multiplicity 1. It remains to show the claim. For this write
 \be
Q_{(\kappa)}^\star&=&V_{(\kappa)}^\star H_{(\kappa)}^\star=\l(\sum_{y,x}H_{_{(\kappa-1)}}^{y}\otimes h_{y,x}\r)\l(\sum_{y',x'}V_{_{(\kappa-1)}}^{y'}\otimes v_{y',x'}\r)\\
&=& \sum_{y,y'}H_{_{(\kappa-1)}}^{y}V_{_{(\kappa-1)}}^{y'} \sum_{x}T_{y,y',x}
\ee
and since this last sum is 1, we have
\beq \label{eq:QP}
Q_{(\kappa)}^\star=\bP_{(\kappa-1)}.
\eq
Also, using \eref{eq:reecriture},
\be
\bP_{(\kappa)}&=&\sum_y Q^y_{(\kappa-1)}\otimes \rho_y
\ee
where $\rho_y=v_{y,\star}h_{y,\star}$. It turns out, that $h_{y,\star}v_{x,\star}=1$ for all $x$ and $y$; hence, any product of the form $v_{y_1,\star}h_{x_1,\star}v_{y_2,\star}h_{y_2,\star}\dots,v_{y_n,\star}h_{y_n,\star}$ equals $v_{y_1,\star}h_{y_n,\star}$. Henceforth, $\kappa\geq 0$, and $m\geq 0$,
\[\l(Q_{(\kappa+2)}^\star\r)^{m+2}=\bP_{(\kappa+1)}^{m+2}=\sum_{x,x'} \l[Q_{\kappa}^x(Q_{\kappa}^\star)^{m}Q_{\kappa}^{x'}\r]\otimes v_{x,\star}h_{x',\star}.\]
Hence, if for some $\kappa\geq 0$, and $m\geq 0$, $(Q_{(\kappa)}^\star)^{m}=(Q_{(\kappa)}^\star)^{m-1}$ then $(Q_{(\kappa+2)}^\star)^{m+2}=(Q_{(\kappa+2)}^\star)^{m+1}.$ The initial conditions  
being $Q_{(0)}^\star=\begin{bmatrix} 1\end{bmatrix}$ (and $Q_{(1)}^\star=\begin{bmatrix} 1\end{bmatrix}$)  we get  $(Q_{(0)}^\star)^1=(Q_{(0)}^\star)^0$, and then $(Q_{(2\kappa)}^\star)^{2\kappa+1}=(Q_{(2\kappa)}^\star)^{2\kappa}$ (and $(Q_{(2\kappa+1)}^\star)^{2\kappa+1}=(Q_{(2\kappa+1)}^\star)^{2\kappa}$ as well).~$\Box$ \medskip \medskip

Denote by $L_{(\kappa)}$ and $R_{(\kappa)}$ the left and right eigenvectors of $Q_{(\kappa)}^\star$ associated with the eigenvalue 1. Since $(Q_{(\kappa)}^\star)^{m}$ converges to $R_{(\kappa)}L_{(\kappa)}$ when $m\to+\infty$, and since $(Q_{(\kappa)}^\star)^{m+1}=(Q_{(\kappa)}^\star)^{m}$, for $m$ large enough,
\beq\label{eq:le}
(Q_{\kappa}^\star)^{m}=R_{(\kappa)}L_{(\kappa)}.
\eq
Moreover  $L_{(\kappa)}$ and $R_{(\kappa)}$ can be normalised such that $L_{(\kappa)}R_{(\kappa)}=1$. Notice in \eref{eq:le} the equality; in similar situations only convergence holds.  For $m$ large enough the quantity of interest
\[\Trace(Q_{(\kappa)}^1(Q_{(\kappa)}^\star)^{m})=\Trace(Q_{(\kappa)}^1R_{(\kappa)}L_{(\kappa)}).\]

Using \eref{eq:QP}
\begin{equation}\label{eq:122}
(Q_{(\kappa)}^\star)^m=\bP_{(\kappa-1)}^m=H_{(\kappa-1)}^\star (Q_{(\kappa-1)}^\star)^{m-1}V_{(\kappa-1)}^\star,
\end{equation}
which leads to
\be
R_{(\kappa)}L_{(\kappa)}&=&V_{(\kappa)}^\star {R}_{{(\kappa+1)}}{L}_{{(\kappa+1)}}H_{(\kappa)}^\star,\\
{R}_{(\kappa)}{L}_{(\kappa)}&=&H_{(\kappa-1)}^\star R_{(\kappa-1)}L_{(\kappa-1)}V_{(\kappa-1)}^\star,
\ee
and then (since all matrices have rank 1 and $L_{(\kappa)}={}^tR_{(\kappa)}$)
\ben\label{eq:link}
{L}_{(\kappa)}&=&H_{(\kappa-1)}^\star L_{(\kappa-1)},\\
              &=&H_{(\kappa-1)}^\star H_{(\kappa-2)}^\star\dots H_{(0)}^\star.
\een
Since $\Trace(Q_{(\kappa)}^1R_{(\kappa)}L_{(\kappa)})=\Trace(H_{(\kappa)}^1R_{(\kappa)} L_{(\kappa)}V_{(\kappa)}^1)$ we have also
\[H_{(\kappa)}^1L_{(\kappa)}=H_{(\kappa)}^1H_{(\kappa-1)}^\star\dots H_{(0)}^\star\]
a triangular product whose computation seems to be quite difficult. We may also note the following
\begin{lem} Let $d^{(n)}(1)=`P(Y^{\Sq(n)}_v=1)$, the density of the gas process of type 2 (this is $\G^{\Sq(n)}(x,y))$ up to change of variables, by Proposition \ref{pro:fo}). For $n$ large enough
\[d^{(n)}(1)= \frac{\sum_{i\geq 0} L_{(\kappa)}[1,2i+1]R_{(\kappa)}[{2i+1,1}]}{\sum_{j\geq 0} L_{(\kappa)}[{1,j}]R_{(\kappa)}[{j,1}]}.\]
\end{lem}
\Proof For short, we don't write the indices ${(\kappa)}$. For $n$ large enough (by \eref{eq:le}) we have $(Q^\star)^{n}= RL$. Hence, 
\[d^{(n)}(1)=\Trace(Q_1(Q^\star)^{n-1})/\Trace((Q^\star)^{n})=\Trace(LV^1 H^1R)/\Trace(LR).\]
Introduce $L(1)=\begin{bmatrix}L[1,i] \1_{i \mod 2=1} \end{bmatrix}$ and $R(1)=\begin{bmatrix}R[{i,1}] \1_{i \mod 2=1} \end{bmatrix}$, the vectors $L$ and $R$, where the even entries are sent to 0. Now, clearly
$LV^1=(LV^\star)(1)$ and  $H^1 R=(H^\star R)(1)$, then $(LV^\star)(1)=L(1)$ and $(H^\star R)(1)=R(1)$.~$\Box$

Let us come back to the matrices $Q_{\infty}^1$ and $Q_{(\infty)}^0$ solution of the rewriting like system \eref{eq:rew1}.
Again, the value of the corner of $Q^x[0,0]$ is given by $\rho_Q[0,0]\bM_{\bq}[0,0]\begin{bmatrix} 1_{x=0}\\ 1_{x=1}\end{bmatrix}$. 
\begin{conj} The matrix $Q^\star_{\infty}$ has eigenvalue 1 with multiplicity 1. Let $L^{\infty}$ and $R^{\infty}$ be the left and right eigenvector, such that $L^{\infty}={}^{t}R^{\infty}$. We have
\begin{equation}\label{eq:inf}
`P(Y_x=1)=\sum_i L^{\infty}[1,2i+1]R^{\infty}[1,2i+1]/ \sum_j L^{\infty}[{1,j}]R^{\infty}[{j,1}].
\end{equation}
 \end{conj}
We have seen that $Q_{(\kappa)}^x\sous{\longrightarrow}{\kappa\to +\infty} Q_{(\infty)}^x$ (simply) and $Q_{(\kappa)}^\star$ has a unique eigenvalue 1, the other ones being 0. This convergence is not sufficient to deduce that the infinite matrix $Q^\star_{\infty}$ has eigenvalue 1, and even if it is the case, the convergence of the numerator and denominator in \eref{eq:inf} may not converge.  In the general case, for 4 given matrices $(q_{x,y}, (x,y)\in\{0,1\}^2)$, the same question arises2: can we find the eigenvectors (and eigenvalues) of the matrix $Q^\star_{\infty}$ that solves the rewriting systems  \eref{eq:rew1}.

\begin{Ques} Is it possible to find solutions to \eref{eq:finite} with infinite matrices $V^1, H^1, V^0, H^0$, in the case $T=T^\bY$ and such that moreover any product of the form $\prod_{i=1}^n V^{x_i}H^{x_i}$ converges ?
\end{Ques}

\section{Other similar but different considerations} 

We present some alternatives to Lemma \ref{pro:dec}. Even if morally what is done has 
more of less the same taste as this Lemma, we were not able to reduce the following considerations to it.

\subsection{Research of a solution on the zigzag}
\label{sec:zz}
We discussed above the construction of an invariant measure relying on the research of a measure of the type $\Trace(\prod_{i=1}^k Q^{x_i})$ on the rows of the cylinder. Two closely related constructions can be proposed. The first one is quite close to that discussed in \cite{BM1} around the question of ``Markovian Field''. 
\begin{figure}[ht]
\psfrag{l0}{$d_0$}
\psfrag{l}{$l$}
\psfrag{r}{$r$}
\psfrag{x_1}{$x_1$}\psfrag{x_2}{$x_2$}\psfrag{x_n}{$x_n$}
\psfrag{y_1}{$y_1$}
\psfrag{y_n}{$y_n$}
\centerline{\includegraphics[height=2cm]{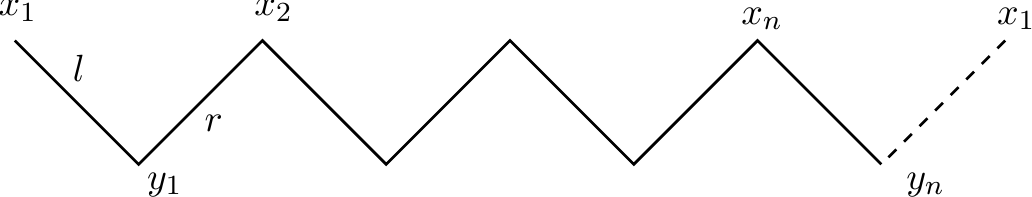}}
\captionn{\label{fig:zig}Decomposition of the measure along the zig-zag.}
\end{figure} The idea is to search for an invariant measure on a zigzag on the cylinder having the following forms (idea used with success by Dhar): for any $(\bx,\by)\in E_n^2$,
\beq\label{eq:formprod}
\mu(\bx,\by)= c_n \prod_{i=1}^n d_{x_i,y_i}u_{y_i,x_{i\oplus1}},
\eq
for some complex numbers $(d_{a,b},u_{a,b},a,b \in\{0,1\})$, where ``$d$'' is chosen for ``down'', and $u$ for ``up''.  But such a measure has a cyclic Markovian structure on the lines since for $m_{a,b}=\sum_c d_{a,c}u_{c,b}$ and $\tilde m_{a,b}=\sum_c u_{a,c}d_{c,b}$ the induced law on the line above (up) is $\mu^{u}$ when the measure below (down) is $\mu^d$ with $\mu^u(\bx)=c_n\prod m_{x_i,x_{i \oplus 1}}$ and $\mu^d(\by)=c_n\prod \tilde{m}_{y_i,y_{i \oplus 1}}$.
 Since $`P(Y_i=y_i | X_i=x_i,X_{i\oplus 1}=x_{i\oplus 1})=T_{x_i,x_{i\oplus 1},y_i}$  a sufficient condition for $\mu^u=\mu^d$ and $\mu^d= \mu^u \bT$ is that for any $a,b,c,x,y$ 
\ben\label{eq:two-equ}
\left\{\begin{array}{l}
{d_{a,c}u_{c,b}}/{m_{a,b}}=T_{a,b,c}, \\
\tilde m_{x,y}= m_{x,y}.\end{array}\right.
\een 
Notice that condition \eref{eq:two-equ} can be weaken a bit since $\mu^u=\mu^d$ does not imply $\tilde m_{x,y} = m_{x,y}$, but rather $\tilde m_{a,b}=w_{a,b} {m}_{a,b}$ for $w_{a,b}$ such that, for any $\bx \in E_n$ (or sufficiently in the support of $\mu^a$), 
\beq\label{eq:prod1}
\prod_{i=1}^n w_{x_i,x_{i\oplus1}}=1.
\eq
For measures with support $E_n$, letting $N_{a,b}(\bx)=\#\{i\in\{0,\dots,n-1\}: (x_i,x_{i\oplus1})=(a,b)\}$, \eref{eq:prod1} rewrites
\[\prod_{(a,b)\in \{0,1\}^2}w_{a,b}^{N_{a,b}(\bx)}=1.\]
Since on $E_n$, $N_{1,0}=N_{0,1}$ and $N_{1,0}+N_{1,1}=\#\{i:x_i=1\}=n-\#\{i:x_i=0\}=n-N_{0,1}-N_{0,0}$, and$N_{1,1}=n-2N_{1,0}-N_{0,0}$, this condition is fulfilled if $w_{0,0}=w_{1,1}=1$, $w_{0,1}w_{1,0}=1$, giving us one degree of freedom. If we deal with measures on $\{0,1,2,\dots, \kappa\}$, it suffices that the product on all finite cycles $w_{x_1,x_2}\dots w_{x_{l-1}x_1}$ equal 1 for $l\leq k$. This provides also some degrees of freedom.
% (for example for $k=2$ this gives  
%\[\left\{ w_{{i,i}}=1 \textrm{ for }i \in\{0,1,2\}, w_{i,j}=1/w_{j,i} \textrm{ for }i< j,
%w_{i,j}w_{j,k}w_{k,i}=1 \textrm{ for } \{i,j,k\}=\{0,1,2\}\r\},\]
% so two degrees of freedom, and 3 if $k=3$).

Hence, the existence of a product form as \eref{eq:formprod} is equivalent to the existence of solutions to  the following system of equations
\ben
\label{eq:two-equ2}
\left\{\begin{array}{l}
\frac{d_{a,c}u_{c,b}}{m_{a,b}}=T_{a,b,c}, \\
 m_{x,y}w_{x,y}= m_{x,y}.\\
w_{0,0}=w_{1,1}=1, w_{0,1}w_{1,0}=1\end{array}\right.
\een
This is a finite algebraic system, and solutions can be found using the computation of a Gröbner basis. Again,  in order to avoid trivial solutions, an additional equation has to be added:  letting 
$M:=(m_{x,y})_{(x,y)\in\{0,1\}^2}$, the equation $\Trace(M^n)=\mu(E_n)$ let us sees that the existence of a non zero eigenvalue for $M$ is necessary and sufficient for non triviality:
\beq\label{eref:CNT}
CNT_2:~~\textrm{ 1 is an eigenvalue of }M.
\eq
The computation of the Gröbner basis of the set of polynomials corresponding to the equations \eref{eq:two-equ2} gives all transitions $T$ for which solution exists (of course, as usual $T^\bY$ is not in this set). \\
One can go further searching for a matrix type solution. It suffices to mix \eref{eq:two-equ} with the matrix considerations of the previous subsection. We then search matrices $\l(D_{a,b},U_{a,b},\textrm{ for }(a,b)\in \{0,1\}^2\r)$ of size say $k\times k$ such that for $m_{a,b}=\sum_c D_{a,c}U_{c,b}$ and  $\tilde m_{a,b}=\sum_c U_{a,c}D_{c,b}$, we have
\ben
\label{eq:two-equ-Mat}
\left\{\begin{array}{l}
D_{a,c}U_{c,b}=T_{a,b,c}m_{a,b}, \\
P m_{x,y}= w_{x,y}\tilde m_{x,y} P 
\end{array}\right.
\een
for some $k\times k$ invertible matrix $P$, $(w_{x,y},x,y\in\{0,1\})$ solution of \eref{eq:prod1}. Again some non-degeneracy conditions must be added: define  $K:=(m_{x,y})_{(x,y)\in\{0,1\}}$ (defined by block, and having size $2k \times 2k$). The needed condition is $\Trace(M^n)\neq 0$, for $n$ large. By triangulation of $M$, it appears clearly that a necessary and sufficient condition is that $M$ has a non zero eigenvalue, which amounts to imposing $CNT_2$.

We were unable to show existence/non-existence of such matrices $(P,D,U)$ and weights $(w_{x,y})$ solution of \eref{eq:two-equ-Mat} + $CNT_2$ in the case $T=T^\bY$ for $m\geq 3$.

\subsection{Research of a solution by projection}

Recall a simple fact~: if $(X_i,i\geq 0)$ is a Markov chain with state spaces $S$ (with $\#S>2$), and if $\phi:S\to S'$, then in general $(\phi(Z_i),i\geq 0)$ is not a Markov chain. Here, one may search if the process $\bY$ (on $E_n$), gas of type 2, whose measure $\mu$ is solution of $\mu^\bY= \mu^\bY\bT^\bY$, can be written $Y_i=\phi(Z_i)$ for some function $\phi$, and some process $(Z_i,1\leq i\leq n)$ taking its values in a state space $S$ with $|S|>2$, and having a simple representation. We were not able to find such a solution.

\section{Extensions}
\label{sec:extensions}
In this section we discuss some extensions of the method we developed above: first to the triangular lattice, second to processes with more than 2 states.

\subsection{Triangular lattice}
\label{sec:tri-lat}
We proceed as in Section \ref{sec:zz} where a zigzag is considered.
Since no new result are provided for DA on the triangular lattice, we just explain how the previous method can be adapted here. First the two probabilistic local transitions needed for the definition of the gases of type 1 and 2 on the triangular lattice are 
\be
T^\bX_{abc,d}&=& `P(B_p(1-a)(1-b)(1-c)=d)\\
T^\bY_{abc,d}&=& `P(B_pabc+(1-B_p)B_q=d).
\ee
Again Proposition \ref{pro:fo} says that the density of the corresponding gas provides up to a change of variables $\G(x,1)$ and $\G(x,y)$, the area, and area-perimeter GF of DA on this lattice.
\begin{figure}[ht]
\psfrag{u_1}{$u_1$}
\psfrag{a}{$a$}
\psfrag{b}{$b$}
\psfrag{c}{$c$}
\psfrag{d}{$d$}
\psfrag{d_1}{$d_1$}
\psfrag{u_2}{$u_2$}\psfrag{d_2}{$d_2$}
\centerline{\includegraphics[height=4cm]{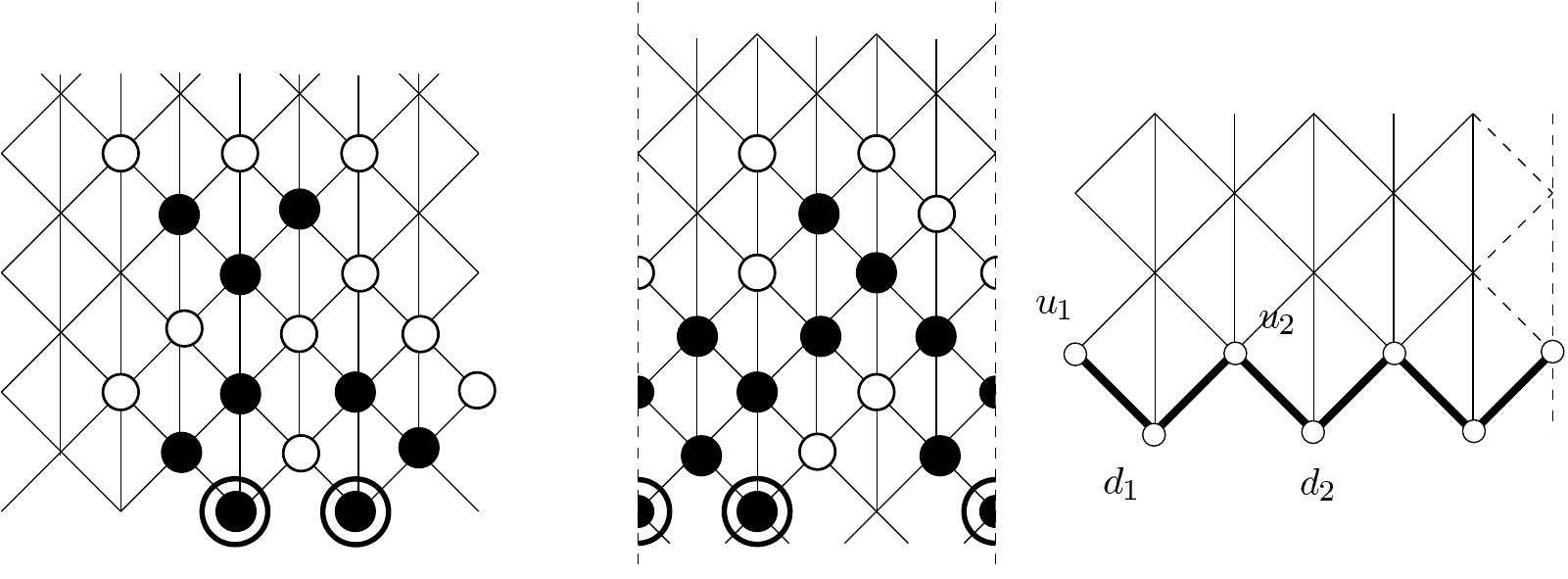}}
\captionn{\label{fig:tri}Representation of the transition. DA on the triangular lattice; the white circles represent perimeter sites. On the second picture, a DA on the triangular lattice with periodic conditions; on the last picture, a zig-zag on the triangular lattice} 
\end{figure}
As explained in Section \ref{sec:zz}, one searches for a representation of the zigzag process distribution as follows:
\beq\label{eq:dec-mes-zz}
\mu_{u_1,d_1,u_2,\dots,u_n,d_n}=\Trace\l(D^{u_1d_1}\dots U^{d_nu_1}\r)
\eq
for some matrices $D^{u,d}$ and $U^{d,u}$ for $u,d\in\{0,1\}^2$ (see Figure \ref{fig:tri} to see the respective  positions of the $u_i$'s and $d_i$'s). To find some finite matrices doing the job, it it sufficient to find matrices $D^{u,d}$ and $U^{d,u}$ for $u,d\in\{0,1\}^2$ solving the following system:
\beq\label{eq:tri-fi}
\left\{
\begin{array}{l}
D^{ab} U^{\bar{b}c}=0,U^{ab}D^{\bar{b}c}=0, \textrm{ for }a,b,c \in \{0,1\}, \\
D^{df} U^{fd'}= \sum_u U^{du}D^{ud'}T_{dud',f}   \textrm{ for } d,f,d' \in \{0,1\}.
\end{array}
\right.
\eq
For $T=T^\bX$ there is a solution with ``matrices'' $D^{x,y},U^{x,y}$ of size 2, for example $D^{x,y}=U^{x,y}$ for any $x,y$, $D^{1,1}=\matz, D^{1,0}=\begin{bmatrix} -1 & -r \\ -(rp+1+2p)/p & 1 \end{bmatrix}$, $ D^{0,1}=\begin{bmatrix} -rp & rp \\ -rp & rp \end{bmatrix}$, $D^{0,0}=\begin{bmatrix} 1 & r \\ 1 & r \end{bmatrix}$, where $r$ satisfies $p+(1-2p)r+r^2p=0$. This leads to the searched density (similar approach are present in \cite{Al,BM1}). \par
Again no such chance arises for $T^\bY$. To adapt the construction of growing matrices as explained in Corollary \ref{cor:big-construction}, 
some single line and single column matrices $h_{a,b,c}, v_{a,b,c}$ indexed by $(a,b,c)\in\{0,1\}^3$ which solves the following system must be found:
\beq\label{eq:tri-lat}
h_{dxz}v_{xd'z'}=T_{dxd',z} 1_{z=z'}.
\eq 
The idea then is to grow the matrices $D_{(\kappa)},U_{(\kappa)}$ as follows:
\[\l\{\begin{array}{l}
D^{dz}_{(\kappa+1)}=\sum_{x}U_{(\kappa)}^{dx}\otimes h_{dxz}\\
U^{zd'}_{(\kappa+1)}=\sum_{x'}D_{(\kappa)}^{x'd'}\otimes v_{x'd'z}
\end{array}\r.\]
Under this condition if $\mu_{(\kappa)}$ has a representation as that given in \eref{eq:dec-mes-zz} with some matrices $D_{(\kappa)},U_{(\kappa)}$ (instead of $D, U$) such that $D_{(\kappa)}^{xy}U_{(\kappa)}^{y',z}=U_{(\kappa)}^{xy}D_{(\kappa)}^{y',z}=\matz$ if $y\neq y'$, then this property is inherited for $D_{(\kappa+1)}, U_{(\kappa+1)}$.
On the zigzag below (as drawn on Figure \ref{fig:tri}), then measure $\mu_{(\kappa+1)}$ will also be given  as in \eref{eq:dec-mes-zz}, with $D_{(\kappa+1)},U_{(\kappa+1)}$ instead of $D_{(\kappa)},U_{(\kappa)}$. Again, solutions to the system \eref{eq:tri-lat} exist: this permits one to grow some matrices $D_{(\kappa)},U_{(\kappa)}$, and then to have a representation of the measure $\mu_{(\kappa)}$. Infinite matrices $U_{(\infty)}^{ab}$ and $D_{(\infty)}^{ab}$ appear again by some passage to the limit. We were not able to deduce from them $\G$. 

\subsection{Bond percolation}

Let $A$ be a DA on $G=(V,E)$. A bond in $A$ is an edge $e\in E$ between elements of $A$; let  $N(A)$ be the number of  bonds in $A$. The GF $B_{C}^G(x,y)=\sum_{A} x^{|A|} y^{\#N(A)}$ of DA with source $C$ counted according to the area and number of bonds can be obtained also by computing the density of some gas process (this is explained in \cite{BM1}, page 21, in the case $G=\Sq(n)$).  Indeed, for any $G\in{\cal G}$,
\beq\label{eq:bp}
B_C^G=x^{|C|}\sum_{D\subset  \Ch(C)}B_{D}^G y^{C\to D}
\eq
where $C\to D$ is the number of bonds between $C$ and $D$. For any cell $d\in \Ch(C)$, let $S_C(d)=\{(c,d)\in E~| c\in C\}$ be the set of bonds from $C$ to $d$. Further, denote by 
$D_{C}(i)=\{d \in D~|~|S_C(d)|= i\}$ be the subset of $D$ of cells being extremities of $i$ edges coming from $C$. We have
\beq\label{eq:bp-2}
B_C^G=x^{|C|}\sum_{D\subset \Ch(C)} B_{D}^G y^{\sum i |D_C(i)|}.
\eq
We will define a gas whose density will coincide up to a change of variables to $B_C$. For this, associate with the set of vertices $V$ i.i.d. Bernoulli$(p)$ random variables (denoted $(B_p^x,x\in V)$), and with the edges of $E$,  i.i.d. Bernoulli$(q)$ random variables (denoted $(B_q^x,x \in E)$.  Consider now the gas defined by 
\beq \label{eq:tran}
X_x= B_p^x\prod_{d\in \Ch(x)}  1- \l(X_d  \prod_{(x,d)\in E}B_q^{(x,d)}\r).
\eq
Taking the expectation in the previous line, leads to 
\be
`P(X_x=1, x\in C)%&=&p^{|C|}`E\l(\prod_{d\in \Ch(C)}  1- \l(X_d  \prod_{(x,d)\in E} B_q^{(x,d)}\r)\r)\\
                 &=&p^{|C|}\sum_{D\subset \Ch(C)} `P(X_x=1,x\in D)(-1)^{|D|} q^{\sum i |D_C(i)|}.
\ee 
Set now $H_C(p,q)=(-1)^{|C|}`P(X_x=1, x\in C)$. The last equation rewrites
\[H_C(p,q)=(-p)^{|C|}\sum_{D\subset \Ch(C)} H_D(p,q) q^{\sum i |D_C(i)|}\]
in other words, the family $(H_C(-p,q),C)$  satisfies the same equations as the family $(B_C(p,q),C)$, the  initial conditions being $B_{\varnothing}=1$ and $H_{\varnothing}=1$. On the square lattice, set
\[ T^\bB_{x,y,z}=`P(X=z | X_1=x,X_2=y)=`P(z= B_p ( 1- x  B_q^{(1)}) ( 1- y  B_q^{(2)}).\] Again, there exist solutions to the equations $h_{x,y}v_{x',y'}=1_{y=y'}T^\bB_{x,x',y}$, and we may also impose that $M_\bq(0,0)^l$ converges. This permits again to use Corollary \ref{cor:big-construction} to represent $\mu_{(\kappa)}$.

\subsection{Bicolouration}
\label{seq:bic}
We present a new model of gas $X$ taking three values 1, 2 or 3 (we see the value as a colour), having an interest from a combinatorial point of view, and illustrating the universality of the present approach. 
Let $G=(V,E)$ be in ${\cal G}$, and let $(C_x,x\in V)$ be a i.i.d. random colouring of the vertices of $G$, such that
$`P(C=i)=p_i$ for $i\leq 2$, and $\sum p_i=1$.
We set
\begin{equation}\label{eq:3coul}
X_x:= C_x \prod_{c\in \Ch(x)}  1_{X_c\neq C_x}.
\end{equation}
The arguments given in Lemma \ref{lem:welldefi} implies that $X$ is a.s. well defined for $p_1+p_2$ small enough.
According to \eref{eq:3coul}, if $C_x=0$ then $X_x=0$ (with proba. 1), and if $C_x=i$ for $i\in\{1,2\}$, then $X_x=i$ with proba. $p_i$ if no child of $x$ has colour $i$, and $X_x=0$ in the other cases.\par

Clearly $`P(X_1=1)=-\G^{G}(-p_1,1)$ and $`P(X=2)=-\G^G(-p_2,1)$ since the gases $Y=1_{X=1}X$ and $Y'=1_{X=2}X/2$ have the same transitions as the gas of type 1 with parameters $p_1$ and $p_2$ (and by Lemma \ref{pro:fo}). 
\begin{figure}[ht]
\centerline{\includegraphics[height=3.5cm]{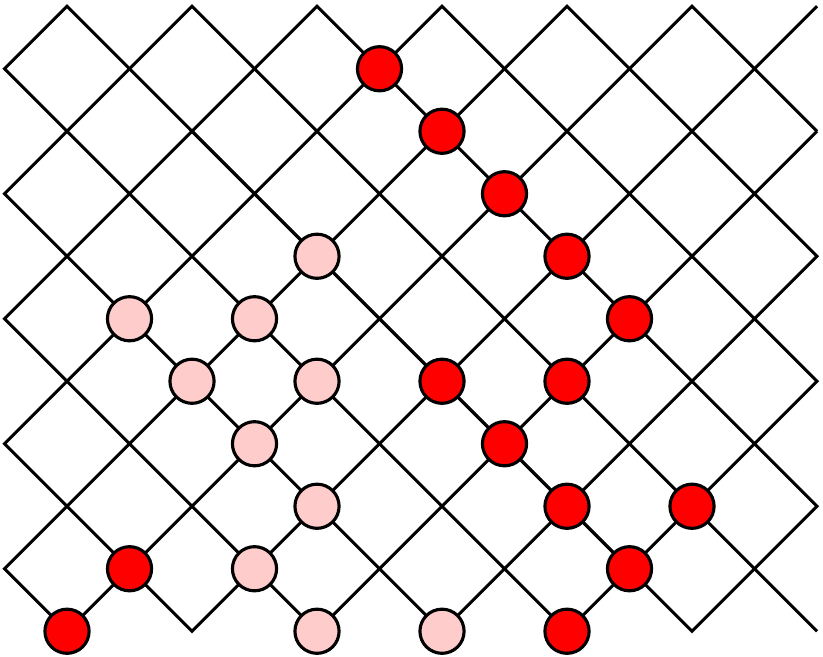}}
\end{figure}
This gas is related to the counting of bi-coloured DA, having no bicoulor neighbouring sites (this model is interesting only when several sources are involved). \par
Formally consider two sets $S_1:=\{c_1,\dots,c_k\}$ and $S_2:=\{d_1,\dots,d_l\}$ such that $S_1\cap S_2=\varnothing$ and such that $S_1\cup S_2$ is a free set. We call bicoloured DA a pair $(A,l)$, where $l:A\to\{1,2\}$ (where $l(a)$ is seen as the colour of $a$). A bicoloured DA is said to be well coloured, if for any $(a,b)\in A^2\cap E$, $l(a)=l(b)$, meaning that neighbouring sites have the same colour.

Denote by  ${\cal A}_{S_1,S_2}$ the set of well coloured DA $(A,l)$ with source $S_1 \cup S_2$ such that $l(S_i)=\{i\}$ (the cells of the sources $S_i$ have colour $i$). This model is hard to  deal with using heap of pieces arguments. Let $\G_{S_1,S_2}$ be the corresponding GF, counted well coloured DA according to their number of cells of each colour. 
\begin{pro} We have
\beq\label{eq:bi-col}
\G_{S_1,S_2}(-p_1,-p_2)=(-1)^{\#S_1+\#S_2}`P(X_x=1, x\in S_1, X_x=2,x\in S_2).
\eq
\end{pro}
Here again, the gas transition $T(a,b,c)=`P(C_x = c ~| C_{c_1}=a, C_{c_2}=b)$ can be written 
$T(a,b,c)=h_{a,c}v_{b,c'}1_{c=c'}$ for some monoline and monocolumn $h$ and $v$, for any $a,b,c,c'\in\{0,1,2\}$. Again, this ensures the existence of a representation of the measure on the $\kappa$th line of the lattice using some matrices $V^x_{(\kappa)},H^x_{(\kappa)},Q^x_{(\kappa)}$, for $x\in\{0,1,2\}$, starting for some measure $\mu_{(0)}(\bx)=\Trace(\prod_{i=0}^n V^{x_i}_{(0)}H^{x_i}_{(0)})$ and some matrices $V^x_{(0)},H^x_{(0)}$, such that $V^x_{(0)}H^y_{(0)}=\matz$ for $x\neq y$.
\par
Here the case is particularly interesting: $`P(X_x=j)$ for $j\in\{0,1,2\}$ as well as $`P(X_x=j, x\in C)$ are easy to compute: the reason is that the projection $Y$ and $Y'$ (as defined above) are well known, and have a simple product form: they correspond in the first case to identify the states $2$ and $0$ and in the second one to identify $0$ and $1$.  Hence, both $Y$ and $Y'$ are Markovian on a line of the lattice (and hard particle model on the zigzag), but  $X$ is not Markovian on the lines (or on the zigzag).  The combinatorial issue is not to find the density of $X$, but rather to compute a quantity as $`P(X_0=1,X_1=2)$. At this moment, I am not able to do this. \medskip\medskip

\noindent \bf Thanks : \rm  Grateful thanks to David Renault for numerous stimulating discussions about this work.  Many thanks are due to the anonymous referees for their comments.

\small
\bibliographystyle{abbrv}

%\newpage
%\tableofcontents

\end{document}